\documentclass[a4paper]{article}

\pdfoutput=1

\usepackage[english]{babel}
\usepackage[utf8]{inputenc}
\usepackage[colorinlistoftodos]{todonotes}
\usepackage{amsmath}
\usepackage{amssymb}
\usepackage{amsthm} 
\usepackage{mathtools}
\usepackage{mathrsfs}
\usepackage{siunitx}
\usepackage{algorithm}
\usepackage{algpseudocode}

\usepackage{subfig}

\usepackage{tikz}

\usepackage{tabularx}
\usepackage{booktabs}

\usepackage{empheq}

\usepackage[toc,page]{appendix}

\usepackage{fixltx2e}
\usepackage{hyperref}
\usepackage{cleveref}

\everymath{\displaystyle}

\newcommand{\VLA}{V_{\text{LA}}}
\newcommand{\VLV}{V_{\text{LV}}}
\newcommand{\VRA}{V_{\text{RA}}}
\newcommand{\VRV}{V_{\text{RV}}}
\newcommand{\VnLA}{V_{\text{0,LA}}}

\newcommand{\VnRA}{V_{\text{0,RA}}}
\newcommand{\VnRV}{V_{\text{0,RV}}}
\newcommand{\PLA}{p_{\text{LA}}}
\newcommand{\PLV}{p_{\text{LV}}}
\newcommand{\PRA}{p_{\text{RA}}}
\newcommand{\PRV}{p_{\text{RV}}}

\newcommand{\EpLA}{E_{\text{LA}}^{\text{pass}}}

\newcommand{\EpRA}{E_{\text{RA}}^{\text{pass}}}
\newcommand{\EpRV}{E_{\text{RV}}^{\text{pass}}}

\newcommand{\EaMaxLA}{E_{\mathrm{LA}}^{\mathrm{act,max}}}

\newcommand{\EaMaxRA}{E_{\mathrm{RA}}^{\mathrm{act,max}}}
\newcommand{\EaMaxRV}{E_{\mathrm{RV}}^{\mathrm{act,max}}}
\newcommand{\QarSYS}{Q_{\text{AR}}^{\text{SYS}}}
\newcommand{\QarPUL}{Q_{\text{AR}}^{\text{PUL}}}
\newcommand{\QvnSYS}{Q_{\text{VEN}}^{\text{SYS}}}
\newcommand{\QvnPUL}{Q_{\text{VEN}}^{\text{PUL}}}
\newcommand{\CarSYS}{C_{\text{AR}}^{\text{SYS}}}
\newcommand{\CarPUL}{C_{\text{AR}}^{\text{PUL}}}
\newcommand{\CvnSYS}{C_{\text{VEN}}^{\text{SYS}}}
\newcommand{\CvnPUL}{C_{\text{VEN}}^{\text{PUL}}}
\newcommand{\ParSYS}{p_{\text{AR}}^{\text{SYS}}}
\newcommand{\ParPUL}{p_{\text{AR}}^{\text{PUL}}}
\newcommand{\PvnSYS}{p_{\text{VEN}}^{\text{SYS}}}
\newcommand{\PvnPUL}{p_{\text{VEN}}^{\text{PUL}}}
\newcommand{\RarSYS}{R_{\text{AR}}^{\text{SYS}}}
\newcommand{\RarPUL}{R_{\text{AR}}^{\text{PUL}}}
\newcommand{\RvnSYS}{R_{\text{VEN}}^{\text{SYS}}}
\newcommand{\RvnPUL}{R_{\text{VEN}}^{\text{PUL}}}
\newcommand{\LarSYS}{L_{\text{AR}}^{\text{SYS}}}
\newcommand{\LarPUL}{L_{\text{AR}}^{\text{PUL}}}
\newcommand{\LvnSYS}{L_{\text{VEN}}^{\text{SYS}}}
\newcommand{\LvnPUL}{L_{\text{VEN}}^{\text{PUL}}}

\newcommand{\QAV}{Q_{\text{AV}}}
\newcommand{\QMV}{Q_{\text{MV}}}
\newcommand{\QTV}{Q_{\text{TV}}}
\newcommand{\QPV}{Q_{\text{PV}}}

\newcommand{\Rmin}{R_{\mathrm{min}}}
\newcommand{\Rmax}{R_{\mathrm{max}}}

\newcommand{\CircOne}{\boldsymbol{c}_1}
\newcommand{\CircTwo}{\boldsymbol{c}_2}

\newcommand{\CircOneInit}{\boldsymbol{c}_{1,0}}

\newcommand{\CircRhs}{\boldsymbol{D}}
\newcommand{\CircSteady}{\boldsymbol{W}}

\newcommand{\fZero}{{\mathbf{f}_0}}
\newcommand{\sZero}{{\mathbf{s}_0}}
\newcommand{\nZero}{{\mathbf{n}_0}}

\newcommand{\EPchim}{\chi_\text{m}}
\newcommand{\EPCm}{C_\text{m}}
\newcommand{\EPIion}{{\mathcal{I}_{\text{ion}}}}

\newcommand{\EPIapp}{{\mathcal{I}_{\text{app}}}}
\newcommand{\EPIappReduced}{{\widetilde{\mathcal{I}}_{\mathrm{app}}}}

\newcommand{\EPIappReducedMax}{{\widetilde{\mathcal{I}}_{\mathrm{app}}^{\mathrm{max}}}}
\newcommand{\EPIappDuration}{t_{\mathrm{app}}}

\newcommand{\Inv}[1]{{\mathcal{I}_{#1}}}

\newcommand{\IIVf}{\Inv{4f}}

\newcommand{\displ}{\mathbf{d}} 

\newcommand{\BCmecCepiN}{{C_\bot^{\text{epi}}}}
\newcommand{\BCmecKepiN}{{K_\bot^{\text{epi}}}}
\newcommand{\BCmecCepiT}{{C_\parallel^{\text{epi}}}}
\newcommand{\BCmecKepiT}{{K_\parallel^{\text{epi}}}}

\newcommand{\GammaBase}{\Gamma_0^{\text{base}}}
\newcommand{\GammaEpi}{\Gamma_0^{\text{epi}}}
\newcommand{\GammaEndo}{\Gamma_0^{\text{endo}}}

\newcommand{\mecF}{\mathbf{F}}

\newcommand{\tenspiola}{\mathbf{P}}
\newcommand{\identity}{\mathbf{I}}

\newcommand{\mecNref}{{\mathbf{N}}}

\newcommand{\Cbar}{\overline{C}}
\newcommand{\bff}{b_{\text{ff}}}
\newcommand{\bss}{b_{\text{ss}}}
\newcommand{\bnn}{b_{\text{nn}}}
\newcommand{\bfs}{b_{\text{fs}}}

\newcommand{\bsn}{b_{\text{sn}}}
\newcommand{\bfn}{b_{\text{fn}}}

\newcommand{\Eff}{E_{\text{ff}}}
\newcommand{\Ess}{E_{\text{ss}}}
\newcommand{\Enn}{E_{\text{nn}}}
\newcommand{\Efs}{E_{\text{fs}}}
\newcommand{\Esf}{E_{\text{sf}}}
\newcommand{\Esn}{E_{\text{sn}}}
\newcommand{\Ens}{E_{\text{ns}}}
\newcommand{\Efn}{E_{\text{fn}}}
\newcommand{\Enf}{E_{\text{nf}}}
\newcommand{\Eab}{E_{\text{ab}}}
\newcommand{\aZero}{\boldsymbol{a}_\text{0}}
\newcommand{\bZero}{\boldsymbol{b}_\text{0}}

\newcommand{\Cai}{{[\text{Ca}^{2+}]_{\text{i}}}}
\newcommand{\Caizero}{{[\text{Ca}^{2+}]_{\text{i,0}}}}

\newcommand{\gammaf}{\gamma_{f}}
\newcommand{\gammas}{\gamma_{s}}
\newcommand{\gamman}{\gamma_{n}}
\newcommand{\lambdaendo}{\lambda_\text{endo}}
\newcommand{\lambdaepi}{\lambda_\text{epi}}
\newcommand{\kendobar}{\Bar{k}_\text{endo}}
\newcommand{\kepibar}{\Bar{k}_\text{epi}}
\newcommand{\kprimebar}{\Bar{k}'}
\newcommand{\FA}{\boldsymbol{F_{A}}}
\newcommand{\FE}{\boldsymbol{F}_{E}}
\newcommand{\CE}{\boldsymbol{C}_{E}}
\newcommand{\PE}{\boldsymbol{P}_{E}}
\newcommand{\muA}{\mu_{A}}
\newcommand{\RFL}{R_{FL}}

\graphicspath{{pictures/}}

\title{Electromechanical modeling of human ventricles with ischemic cardiomyopathy: numerical simulations in sinus rhythm and under arrhythmia}

\author{Matteo Salvador$^1$,
        Marco Fedele$^1$,
        Pasquale Claudio Africa$^1$,
        Eric Sung$^3$,\\
        Luca Ded\`{e}$^1$,
        Adityo Prakosa$^3$,
        Jonathan Chrispin$^4$,\\
        Natalia Trayanova$^3$,
        Alfio Quarteroni$^{1, 2}$}

\date{$^1$ \footnotesize MOX-Dipartimento di Matematica, Politecnico di Milano, Milan, Italy \\
      $^2$ Professor Emeritus, \'Ecole Polytechnique F\'ed\'erale de Lausanne, Lausanne, Switzerland \\
      $^3$ \footnotesize Department of Biomedical Engineering, Johns Hopkins University, Baltimore, MD, USA, \\
      $^4$ Department of Medicine, Johns Hopkins Hospital, Baltimore, MD, USA \\[2ex]}

\begin{document}
\maketitle

\begin{abstract}
We developed a novel patient-specific computational model for the numerical simulation of ventricular electromechanics in patients with ischemic cardiomyopathy (ICM). This model reproduces the activity both in sinus rhythm (SR) and in ventricular tachycardia (VT).
The presence of scars, grey zones and non-remodeled regions of the myocardium is accounted for by the introduction of a spatially heterogeneous coefficient in the 3D electromechanics model.
This 3D electromechanics model is firstly coupled with a 2-element Windkessel afterload model to fit the pressure-volume (PV) loop of a patient-specific left ventricle (LV) with ICM in SR.
Then, we employ the coupling with a 0D closed-loop circulation model to analyze a VT circuit over multiple heartbeats on the same LV.
We highlight similarities and differences on the solutions obtained by the electrophysiology model and those of the electromechanics model, while considering different scenarios for the circulatory system.
We observe that very different parametrizations of the circulation model induce the same hemodynamical considerations for the patient at hand.
Specifically, we classify this VT as unstable.
We conclude by stressing the importance of combining electrophysiological, mechanical and hemodynamical models to provide relevant clinical indicators in how arrhythmias evolve and can potentially lead to sudden cardiac death.
\end{abstract}


\section{Introduction}
\label{sec: introduction}

Ventricular tachycardia (VT) is a life-threatening arrhythmia that predisposes patients to sudden cardiac death (SCD). Clinically, patients with greater ventricular dysfunction are more likely to develop severe VTs. Even though left ventricular ejection fraction (LVEF) poorly reflects the  mechanical aspect of arrythmias, LVEF$<$35\% has been used as the main metric to determine the level of SCD risk in the clinics \cite{Epstein2008}. Prior work suggests that this greater arrhythmia risk could be due to abnormal electromechanical feedback, which results in greater electrical instability and arrhythmogenic propensity \cite{Taggart1999}. The cardiac mechano-electric coupling presents proarrhythmic effects in pathological scenarios and might induce extra stimuli, early afterdepolarizations or delayed afterdepolarizations \cite{Kohl2011, Timmermann2017AnIA}. However, it is unknown from a clinical perspective how the electromechanical feedback effects manifest during VT, and it is difficult to dissect its mechanisms in a clinical study.

Multiscale, computational heart models incorporating information from the cellular to whole-organ levels are well-suited to dissect fundamental mechanisms of arrhythmogenesis. Biophysically-detailed models of cardiac electrophysiology are well-established and have been utilized widely in the simulation of VTs in both experimental and clinical contexts \cite{Arevalo, Prakosa}. Numerical simulations circumvent experimental and clinical limitations by giving the possibility to test different hypotheses and gain insights on the electromechanical function \cite{Gurev2011}. However, very few studies have been able to create models that synthesize both electrical and mechanical activities into a cohesive framework \cite{Plank2016, Peirlinck2021}. Furthermore, a missing component in the electromechanical studies is the representation of tachycardias and the effect of heart contraction in the re-entrant circuits that underlie arrhythmias.

The goal of this study is to develop a patient-specific computational model of electromechanical activity in the human ventricles with ischemic cardiomyopathy (ICM) and to use this model to study how mechanical contraction affects VT arrhythmia dynamics. This patient-specific electromechanical model would allow us to provide non-invasive, personalized assessments of VT circuits and the corresponding hemodynamical consequences. Electromechanical simulations provide a more accurate knowledge and additional information on VTs with respect to electrophysiological ones. Moreover, their outcomes differ in terms of action potential propagation and conduction velocity once the VT is triggered. By exploiting the coupling with a closed-loop circulation model \cite{Regazzoni2020PartI, Regazzoni2020PartII}, we are also able to differentiate between hemodynamically stable VTs, which can be managed by anti-arrhythmic medication, and hemodynamically unstable ones, that generally require cardioversion. 

\section{Model development}
\label{sec: mathematicalmodeling}
In this section we provide an overview of the different core models involved in the framework of cardiac electromechanics. For more information about the mathematical equations we refer to \cite{GerbiSegregated, GerbiPhD, Gerbi, regazzoni2020thesis, regazzoni2019mor-sarcomeres, Regazzoni2020PartI, Regazzoni2020PartII, Salvador2020}.

\subsection{Electrophysiology}
\label{sec: electrophysiology}

\begin{figure}[t!]
\advance\leftskip-1.2cm
\includegraphics[keepaspectratio, width=1.2\textwidth]{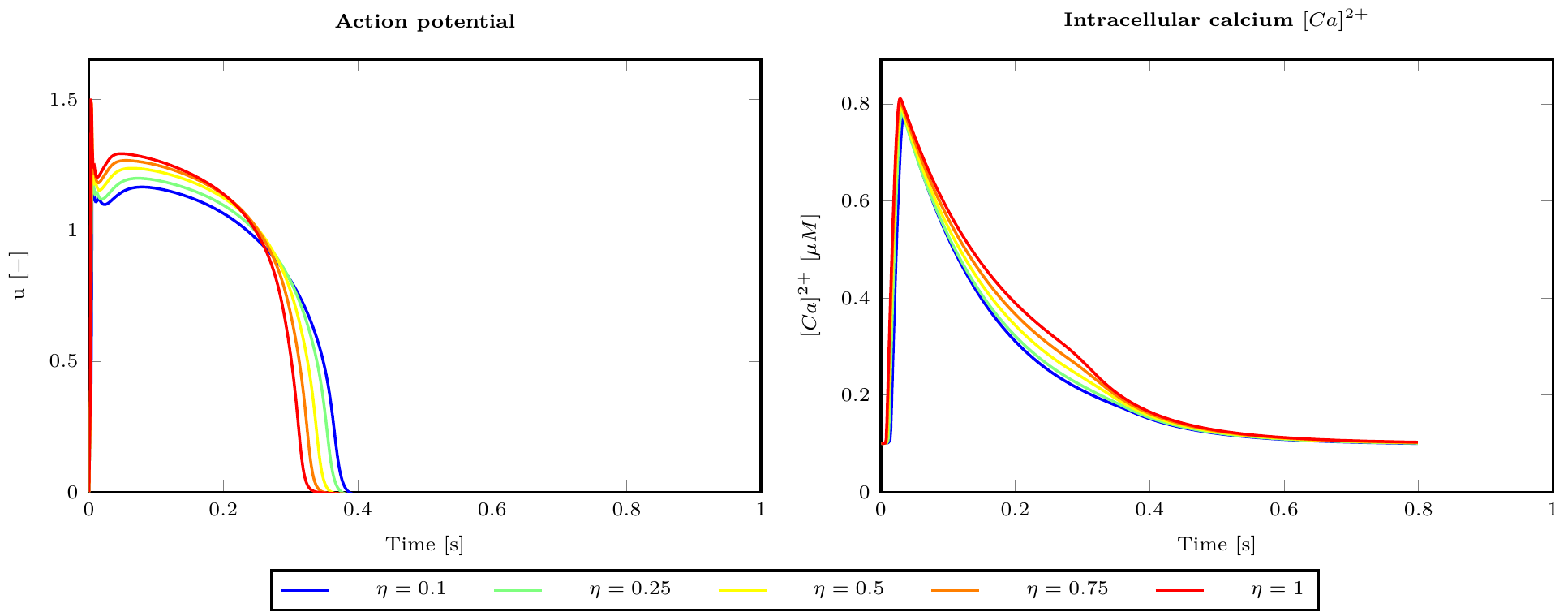}
\caption{Action potential (left) and intracellular calcium concentration $[Ca]^{2+}$ (right) for different values of $\eta=\eta(\boldsymbol{x})$ (TTP06 model, epicardium). $\eta=1$ corresponds to healthy tissue, whereas $\eta \in \{0.1, 0.25, 0.5, 0.75\}$ defines different grey zones. }
\label{fig: AP_Ca_TTP06}
\end{figure}

We model the electrophysiological behavior by means of the monodomain equation, a diffusion-reaction partial differential equation that describes the electric properties of cardiac muscle cells, assuming the same anisotropy ratios between the intracellular and extracellular spaces \cite{Potse, SainteMarie}.
It is a homogenized continuum model, which means that it is used to capture average properties of many cardiomyocytes, and not the behavior of single cells.
We couple it with the ten Tusscher-Panfilov (TTP06) ionic model, due to our focus on the human left ventricle (LV) \cite{TTP06}.
This model permits to describe the microscopic details of the single cardiomyocyte in an accurate and detailed manner \cite{ColliFranzone1}.

The electrophysiology model reads:
\begin{equation}
\begin{cases}
\EPchim \left[ \EPCm \dfrac{\partial u}{\partial t} + \EPIion(u, \boldsymbol{w}, \boldsymbol{z}) \right] = \nabla \cdot ( J \boldsymbol{F}^{-1} \boldsymbol{D}_M \boldsymbol{F}^{-T} \nabla u) + \EPIapp(t) & $in$ \; \Omega_{0} \times (0, T), \\[0.2cm]
\dfrac{\partial \boldsymbol{w}}{\partial t} - \boldsymbol{H}(u, \boldsymbol{w}) = \boldsymbol{0} & $in$ \; \Omega_{0} \times (0, T), \\[0.2cm]
\dfrac{\partial \boldsymbol{z}}{\partial t} - \boldsymbol{G}(u, \boldsymbol{w}, \boldsymbol{z}) = \boldsymbol{0} & $in$ \; \Omega_{0} \times (0, T), \\[0.2cm]
\left( J \boldsymbol{F}^{-1} \boldsymbol{D}_M \boldsymbol{F}^{-T} \nabla u \right) \cdot \boldsymbol{N} = 0 & $on$ \; \partial \Omega_0 \times (0, T), \\
u = u_0 & $in$ \; \Omega_0 \times \{0\}.
\end{cases}
\label{eqn: monodomainionic}
\end{equation}
$\Omega_0 \subset \mathbb{R}^3$ is the computational domain in the reference configuration \cite{Regazzoni2020PartI, Regazzoni2020PartII}: in our application, it is represented by a stress-free patient-specific LV with ICM. $T > 0$ is the final time. $\EPCm$ is the total membrane capacitance and $\EPchim$ is the area of cell membrane per tissue volume. $u$ is the dimensionless transmembrane potential, vector $\boldsymbol{w} = \{w_1, w_2, ..., w_k\}$ expresses $k$ recovery (or gating) variables, which play the role of probability density functions and model the fraction of open ionic channels across the membrane of a single cell, and vector $\boldsymbol{z} = \{z_1, z_2, ..., z_m\}$ defines $m$ concentration variables of specific ionic species (such as calcium $Ca^{2+}$, which plays a major role in heart contraction and mechanical activation). $\EPIapp(t)$ is an external applied current, which simulate in our case the behavior of the Purkinje network \cite{Vergara, Romero}. Indeed we use it to trigger the action potential in specific points of the myocardium. $\EPIion(u, \boldsymbol{w}, \boldsymbol{z})$ is the feedback from the cellular scale into the tissue one, and strictly depends on the chosen ionic model. A Neumann boundary condition is applied all over the boundary and defines the condition of electrically isolated domain. $\boldsymbol{H}(u, \boldsymbol{w})$ and $\boldsymbol{G}(u, \boldsymbol{w}, \boldsymbol{z})$ keep into account the specific features of the TTP06 ionic model \cite{TTP06}. $\boldsymbol{D}_M = \eta \sigma_t \boldsymbol{I} + \eta (\sigma_l - \sigma_t) \fZero \otimes \fZero + \eta (\sigma_n - \sigma_t) \nZero \otimes \nZero$ refers to the diffusion tensor, being $\fZero$ the vector field expressing the fibers direction, $\nZero$ the vector field that indicates the crossfibers direction, and $\sigma_l, \sigma_t, \sigma_n \in \mathbb{R}^+$ the longitudinal, transversal and normal conductivities, respectively \cite{Saffitz}. The parameter $\eta = \eta(\boldsymbol{x}) $ takes into account the effect of ischemic regions both at the macroscopic scale and at the microscopic one. With respect to the latter point, $G_{\text{Na}}(\boldsymbol{x})$, $G_{\text{CaL}}(\boldsymbol{x})$, $G_{\text{kr}}(\boldsymbol{x})$ and $G_{\text{ks}}(\boldsymbol{x})$ conductances of the TTP06 ionic model \cite{TTP06} vary in space according to the following laws:
\begin{equation*}
G_{\text{Na}}(\boldsymbol{x}) = \left[0.38 + \dfrac{10}{9} (\eta(\boldsymbol{x}) - 0.1)(1 - 0.38)\right] G_{\text{Na}},
\end{equation*}
\begin{equation*}
G_{\text{CaL}}(\boldsymbol{x}) = \left[0.31 + \dfrac{10}{9} (\eta(\boldsymbol{x}) - 0.1)(1 - 0.31)\right] G_{\text{CaL}},
\end{equation*}
\begin{equation*}
G_{\text{kr}}(\boldsymbol{x}) = \left[0.30 + \dfrac{10}{9} (\eta(\boldsymbol{x}) - 0.1)(1 - 0.30)\right] G_{\text{kr}},
\end{equation*}
\begin{equation*}
G_{\text{ks}}(\boldsymbol{x}) = \left[0.20 + \dfrac{10}{9} (\eta(\boldsymbol{x}) - 0.1)(1 - 0.20)\right] G_{\text{ks}},
\end{equation*}
where $G_{\text{Na}}$, $G_{\text{CaL}}$, $G_{\text{kr}}$ and $G_{\text{ks}}$ are defined as in \cite{TTP06}. We can potentially consider a continuum of grey zones using linear interpolation, for $\eta(\boldsymbol{x}) \in [0.1, 1]$, going from the full healthy case ($\eta(\boldsymbol{x}) = 1$) to the grey zone described in \cite{Arevalo} ($\eta(\boldsymbol{x}) = 0.1$).
We focus on the modeling of a single possible grey zone, \textit{i.e.} the one reported in \cite{Arevalo}, where $G_{\text{Na}}$, $G_{\text{CaL}}$, $G_{\text{kr}}$ and $G_{\text{ks}}$ conductances are reduced to 38$\%$, 31$\%$, 30$\%$ and 20$\%$ of their physiological values, respectively. This leads to a reduced upstroke and a longer plateau of the action potential. We finally model scars as myocardial regions where no evolution of both the transmembrane potential and all ionic variables occur, where therefore Eq.~\eqref{eqn: monodomainionic} is not actually solved. We report in Fig.~\ref{fig: AP_Ca_TTP06} the evolution over time of the transmembrane potential $u$ and the intracellular calcium concentration $Ca^{2+}$ for different values of the parameter $\eta$.

By defining $\boldsymbol{X}$ and $\boldsymbol{x}$ as the reference and deformed coordinates respectively, we introduce the deformation tensor $\boldsymbol{F} = \boldsymbol{I} + \dfrac{\partial \boldsymbol{d}}{\partial \boldsymbol{X}}$ (with $J=\text{det}(\boldsymbol{F})>0)$. Indeed, in the formulation of our model, we consider the so called mechano-electric feedback (MEF), which models the effect of mechanical deformation $\boldsymbol{F}$ on cardiac electrophysiology \cite{Trayanova}. MEF can be classified in two different categories \cite{Collet2015}. There are both geometrical and physiological contributions. The latter category mainly coincides with stretch-activated channels (SACs) and mechanical modulation of cellular $Ca^{2+}$ handling. In this work, as it can be seen from Eq.~\eqref{eqn: monodomainionic}, we focus on the geometrical part only.

\subsection{Mechanical activation}
\label{sec: mechanicalactivation}
Mechanical activation bridges electrophysiology and passive mechanics.
There are two approaches available in literature, the active stress \cite{Land2017, regazzoni2018redODE} and active strain \cite{Ambrosi2, Ambrosi1} models.
With the former approach, the underlying hypothesis is that an active force is generated by the myocardium, whereas in the latter an active deformation is prescribed to the cardiac tissue.
We propose here modeling of grey zones and scars within the active strain framework \cite{Barbarotta, Rossi1}.

We consider a phenomenological law that keeps into account the local shortening of the fibers $\gammaf$ at the macroscopic level \cite{Azzolin2020, Gerbi, Rossi1, RuizBaier}.
Myocardial displacement $\displ$ and concentration of intracellular calcium ions $\Cai$ play an important role in the time evolution of $\gammaf$.

The phenomenological law reads:
\begin{equation}
\begin{cases}
\dfrac{\partial \gammaf}{\partial t} - \eta \dfrac{\varepsilon}{g(\Cai)} \Delta \gammaf = \eta \dfrac{1}{g(\Cai)} \Phi(\Cai, \gammaf, \displ) & $in$ \; \Omega_0 \times (0, T), \\
\nabla \gammaf \cdot \boldsymbol{N} = 0 & $on$ \; \partial \Omega_0 \times (0, T), \\
\gammaf = 0 & $in$ \; \Omega_0 \times \{0\},
\end{cases}
\label{eqn: activestrain}
\end{equation}
where $g(s) = \muA (\Cai)^2$, $\Phi(\Cai, \gammaf, \displ) = \alpha H_{\text{s}_\text{0}}(\Cai) (\Cai - \Caizero)^2 \\ \RFL(\IIVf) + \sum_{j = 1}^5 (-1)^j (j + 1) (j + 2) \IIVf \gammaf^j$ is the active force and $\RFL(\IIVf)$ is a truncated Fourier series expressing the sarcomere force-length relationship \cite{Gordon}.
Both $\alpha$ and $\muA$ should be calibrated according to the specific case under investigation.
The active deformation is computed by exploiting the following orthotropic form \cite{Ambrosi1}:
\begin{equation}
\FA = \boldsymbol{I} + \gammaf \fZero \otimes \fZero + \gammas \sZero \otimes \sZero + \gamman \nZero \otimes \nZero,
\label{eqn: FA}
\end{equation}
where $\sZero$ and $\nZero$ represent sheets and their normal direction respectively, with $\gammas$ and $\gamman$ corresponding to local shortening or elongation \cite{Barbarotta, Omens}:
\begin{equation}
\gamman = \kprimebar \left( \kendobar \dfrac{\lambda - \lambdaepi}{\lambdaendo - \lambdaepi} + \kepibar \dfrac{\lambda - \lambdaendo}{\lambdaepi - \lambdaendo} \right) \left( \dfrac{1}{\sqrt{1 + \gammaf}} - 1 \right),
\end{equation}
\begin{equation}
\label{eqn: gammas}
\gammas = \dfrac{1}{(1 + \gammaf) (1 + \gamman)} - 1.
\end{equation}
Here $\lambda$ represents a transmural coordinate, varying from $\lambdaendo$ at the endocardium to $\lambdaepi$ at the epicardium, which permits to have a transversely non-homogeneous thickening of the left ventricular wall. $\gammas$ set like \eqref{eqn: gammas} yield to $\text{det}(\FA) = 1$. Coefficient $\eta=\eta(\boldsymbol{x})$ keeps into account mathematical modeling of healthy tissue and a continuum of grey zones, according to its value. As described above, we will consider one single grey zone ($\eta(\boldsymbol{x}) = 0.1$), whereas full healthy tissue is modelled using $\eta(\boldsymbol{x}) = 1$. The reduced value of $\eta(\boldsymbol{x})$ for grey zones induces slower activation and deactivation with respect to the healthy case. Moreover, $\eta(\boldsymbol{x}) = 0.1$ prescribes a lower peak value of $\gammaf$ during the cardiac cycle. These outcomes are obtained by reducing the effects of the diffusion term and the forcing term in Eq.~\ref{eqn: activestrain} respectively.

\subsection{Active and passive mechanics}
\label{sec: mechanics}
We use a nearly-incompressible formulation by weakly penalizing large volumetric variations to define the evolution of the displacement in the myocardium \cite{Simo}.
We model the properties of the cardiac tissue by means of fibers $\fZero$, sheets $\sZero$ and their normals $\nZero$, which permit to obtain highly anisotropic internal stresses associated with a prescribed deformation \cite{Guccione1}.

The momentum conservation equation with boundary and initial conditions reads \cite{Regazzoni2020PartI, Regazzoni2020PartII}:
\begin{equation}
\begin{cases}
\rho_{\text{s}} \dfrac{\partial^2 \displ}{\partial t^2} - \nabla \cdot \tenspiola(\displ, \gammaf) = \boldsymbol{0} & $in$ \; \Omega_0 \times (0, T), \\
(\mecNref \otimes \mecNref) \left( \BCmecKepiN \displ+\BCmecCepiN \dfrac{\partial \displ}{\partial t} \right) \\ \; + \; (\identity - \mecNref \otimes \mecNref) \left( \BCmecKepiT \displ+\BCmecCepiT \dfrac{\partial \displ}{\partial t} \right) + \tenspiola(\displ, \gammaf) \mecNref = \boldsymbol{0} & $on$ \; \GammaEpi \times (0, T), \\
\tenspiola(\displ, \gammaf) \mecNref = \displaystyle \frac{| J \mecF^{-T} \mecNref |}{\int_{\GammaBase} | J \mecF^{-T} \mecNref | d \Gamma_0}\int_{\GammaEndo} \PLV(t) J \mecF^{-T} \mecNref d \Gamma_0 & $on$ \; \GammaBase \times (0, T), \\
\tenspiola(\displ, \gammaf) \mecNref = -\PLV(t) J \mecF^{-T} \mecNref & $on$ \; \GammaEndo \times (0, T), \\[0.2cm]
\displ = \displ_0, \, \dfrac{\partial \displ}{\partial t} = \Dot{\displ}_0 & $in$ \; \Omega_0 \times \{0\}.
\end{cases}
\label{eqn: mechanics}
\end{equation}
A Robin boundary condition is prescribed at the epicardium to account for the effect of the pericardial sac, so that the presence of the pericardium is addressed and modelled \cite{Pfaller}. $\BCmecKepiN$, $\BCmecKepiT$, $\BCmecCepiN$, $\BCmecCepiT \in \mathbb{R}^+$ are local values of stiffness and viscosity constants of the epicardial tissue in the normal or tangential directions,  respectively. At the base of the LV we impose an energy consistent boundary condition that models the effect of the blood flow coming from the left atrium \cite{Regazzoni2020PartI, Regazzoni2020PartII}. $\PLV(t)$ is the internal pressure of the ventricular chamber.
The Piola-Kirchhoff strain tensor $\tenspiola = \tenspiola(\displ, \gammaf)$ incorporates both the passive and active mechanical properties of the tissue. After defining the symmetric positive definite right Cauchy-Green tensor $\boldsymbol{C}=\mecF^T \mecF$, being $\mecF = \boldsymbol{I} + \nabla \displ$ the deformation tensor, the strain energy function $\mathcal{W}$ provides a link between the strain and the energy of the material.
Under the hyperelasticity assumption, the strain energy function can be differentiated with respect to the deformation tensor $\boldsymbol{F}$ to obtain $\boldsymbol{P}$:
\begin{equation}
\tenspiola(\displ, \gammaf) = \dfrac{\partial \mathcal{W}(\boldsymbol{C})}{\partial \boldsymbol{F}}.
\label{eq: Piola}
\end{equation}
In the active strain framework, in addition to the reference configuration $\Omega_0$ and the deformed one $\Omega$, we introduce an intermediate state $\hat{\Omega}$, which represents the active part of the deformation \cite{Ambrosi2, Ambrosi1, Nobile, Rossi1}. The $2^{nd}$ order tensor $\FA$ maps $\Omega_0$ into $\hat{\Omega}$, whereas the $\FE$ tensor transforms $\hat{\Omega}$ into $\Omega$. We finally reach the multiplicative decomposition of $\boldsymbol{F}$ = $\FE \FA$.
The first Piola-Kirchhoff strain tensor $\boldsymbol{P}$ reads:
\begin{equation}
\tenspiola = \text{det}(\FA) \PE {\FA}^{-T}, \quad \quad \PE = \dfrac{\partial \mathcal{W}(\CE, J)}{\partial \FE}.
\end{equation}
For additional details on the final form of tensor $\boldsymbol{P}$ in the active strain framework, we refer the reader to \cite{Gerbi}. \\

Several models are available in literature to describe the anisotropic nature of the tissue, such as the Guccione \cite{Guccione2} or the Holzapfel-Ogden laws \cite{Holzapfel}. In this work we employ the Guccione constitutive law, where the energy function reads \cite{Guccione2, Guccione1}:
\begin{equation}
\begin{split}
\mathcal{W}(\boldsymbol{C}) &= \dfrac{\Cbar}{2} \left( e^Q  - 1 \right) \\
Q &= \bff \Eff^2  + \bss \Ess^2 + \bnn \Enn^2 \\
&+ \bfs \left( \Efs^2 + \Esf^2 \right) + \bfn \left( \Efn^2 + \Enf^2 \right) + \bsn \left( \Esn^2 + \Ens^2 \right).
\end{split}
\label{eqn: Guccione_pt1}
\end{equation}
$\Eab = \boldsymbol{E} \aZero \cdot \bZero$ for $a, b \in \{ f, s, n \}$ are the entries of $\textbf{E} = \dfrac{1}{2} \left( \boldsymbol{C} - \boldsymbol{I} \right)$, i.e the Green-Lagrange strain energy tensor. $\Cbar$ is defined as follows:
\begin{equation}
\overline{C} = C [\eta + (1 - \eta) 4.56] \;\;\;\;\; \eta \in [0, 1],
\label{eqn: Guccione_pt2}
\end{equation}
being $C$ a coefficient fitted from experiments \cite{Guccione1}.

\noindent We assign $\eta = 1$ to healthy tissue, $\eta = 0.1$ corresponds to the single type of grey zone that we consider in this work (as in electrophysiology and activation), and $\eta = 0$ models scar regions. In this way, we model a stiffer myocardium for infarcted areas. We introduce a convex term $\mathcal{W}_{\text{vol}}(J) = \dfrac{B}{2} \left( J - 1 \right) \text{log}(J)$ into the energy function $\mathcal{W}$, such that $J=1$ is its global minimum. This term sets a nearly-incompressible constraint \cite{Cheng, Doll, Yin}. $B \in \mathbb{R}^+$ is the bulk modulus, which has a role in the torsion mechanism of the ventricle and enforces the incompressibility constraint \cite{Gerbi}.

\subsection{Blood circulation and Windkessel afterload model}
\label{sec: circulation}

\begin{figure}[t!]
	\centering
	\includegraphics[width=0.85\textwidth]{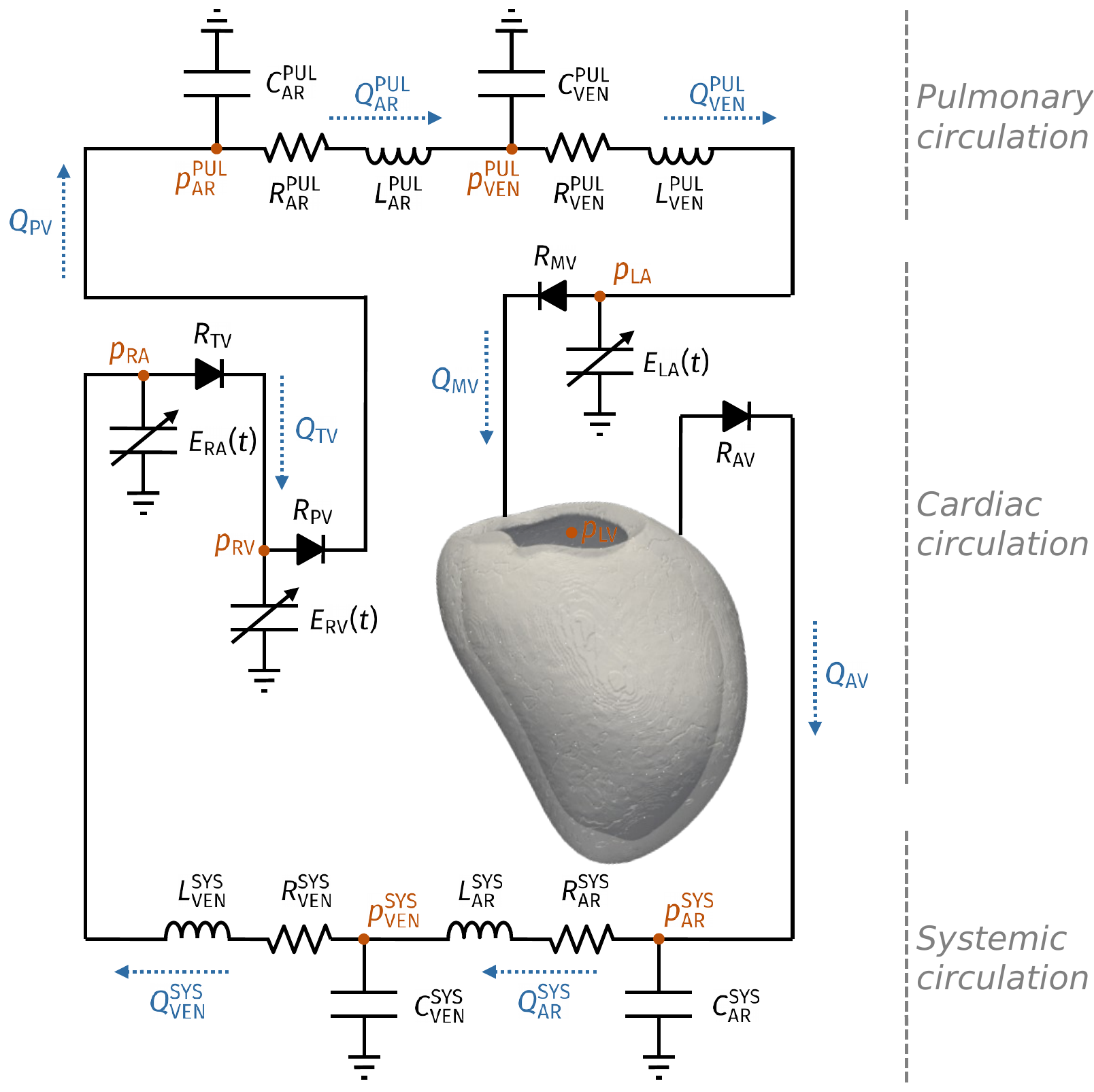}
	\caption{3D-0D coupling between the 3D electromechanics model of a LV with ICM and the 0D closed-loop circulation model. The state variables corresponding to pressures and fluxes are depicted in orange and blue, respectively.}
	\label{fig: 3D0Dcoupling}
\end{figure}

In \cite{Regazzoni2020PartI, Regazzoni2020PartII} we present a 0D closed-loop circulation model of the whole cardiovascular system, previously introduced in \cite{Caruel, Hirschvogel, Washio}. This model is briefly recalled here below in compact form:
\begin{equation} \label{eqn: ode_circ}
\begin{cases}
\dfrac{\mathrm{d} \CircOne(t)}{\mathrm{d} t} = \CircRhs(t, \CircOne(t), \CircTwo(t)) & \qquad t \in (0, T], \\
\CircTwo(t) = \CircSteady(t, \CircOne(t)) & \qquad t \in [0, T], \\
\CircOne(0) = \CircOneInit,
\end{cases}
\end{equation}
with $\CircRhs$ and $\CircSteady$ vector functions defined in \cite{Regazzoni2020PartI, Regazzoni2020PartII}. Variables $\CircOne(t)$ and $\CircTwo(t)$ refer to pressures, volumes and fluxes of the different compartments composing the vascular network:
\begin{equation*}
	\begin{split}
		\CircOne(t) =
		(&\VLA(t), \VLV(t), \VRA(t), \VRV(t), \ParSYS(t), \PvnSYS(t), \ParPUL(t), \PvnPUL(t), \\&\QarSYS(t), \QvnSYS(t), \QarPUL(t), \QvnPUL(t))^T,
		\\
		\CircTwo(t) =
		(&\PLV(t), \PLA(t), \PRV(t), \PRA(t), \QMV(t), \QAV(t), \QTV(t), \QPV(t))^T,
	\end{split}
\end{equation*}
where $\PLA$, $\PRA$, $\PLV$, $\PRV$, $\VLA$, $\VRA$, $\VLV$ and $\VRV$ are pressures and volumes in left atrium, right atrium, left ventricle and right ventricle, respectively; $\ParSYS$, $\QarSYS$, $\PvnSYS$, $\QvnSYS$, $\ParPUL$, $\QarPUL$, $\PvnPUL$ and $\QvnPUL$ define pressures and flow rates of the systemic and pulmonary circulation (arterial and venous), respectively; $\QMV$, $\QAV$, $\QTV$ and $\QPV$ express the flow rates through mitral, aortic, tricuspid and pulmonary valves, respectively. For the complete mathematical description of the lumped 0D circulation model we refer to \cite{Regazzoni2020PartI, Regazzoni2020PartII}.

We couple the 3D electromechanics model of the LV with ICM with the 0D circulation model for the remaining part of the cardiovascular system. We solve the mechanical problem (Eq.~\ref{eqn: mechanics}) with a volumetric constraint on $\VLV(t)$, which is imposed by means of $\PLV(t)$.
For more information about the 3D-0D coupling strategy we refer to \cite{Regazzoni2020PartI, Regazzoni2020PartII}.
In Fig.~\ref{fig: 3D0Dcoupling} we display the electric analog circuit corresponding to the 0D circulation model, along with the coupling with a 3D electromechanical description of a LV with ICM.

In order to ease parameters tuning and exploit the available clinical data on the LV, for the numerical simulations in sinus rhythm (SR) we use a simplified 0D model, \textit{i.e.} a 2-element Windkessel afterload model, which defines the evolution in time of variables $\PLV(t)$ and $\VLV(t)$ only. This model is extensively presented in \cite{GerbiPhD, regazzoni2020thesis, Salvador2020}, along with its numerical discretization. For the sake of completeness, we report here below its equation:
\begin{equation}
C \dfrac{\mathrm{d}\PLV}{\mathrm{d}t} = - \dfrac{\PLV}{R} - \dfrac{\mathrm{d}\VLV}{\mathrm{d}t},
\end{equation}
with $ t \in [0, T]$. $C, R > 0$ are two parameters representing the capacitance and the resistance of the electric circuit that mimics the blood flowing in the aorta.

\section{Numerical discretization}
\label{sec: numericaldiscretization}

\begin{figure}[t!]
\centering
\includegraphics[keepaspectratio, width=1.0\textwidth]{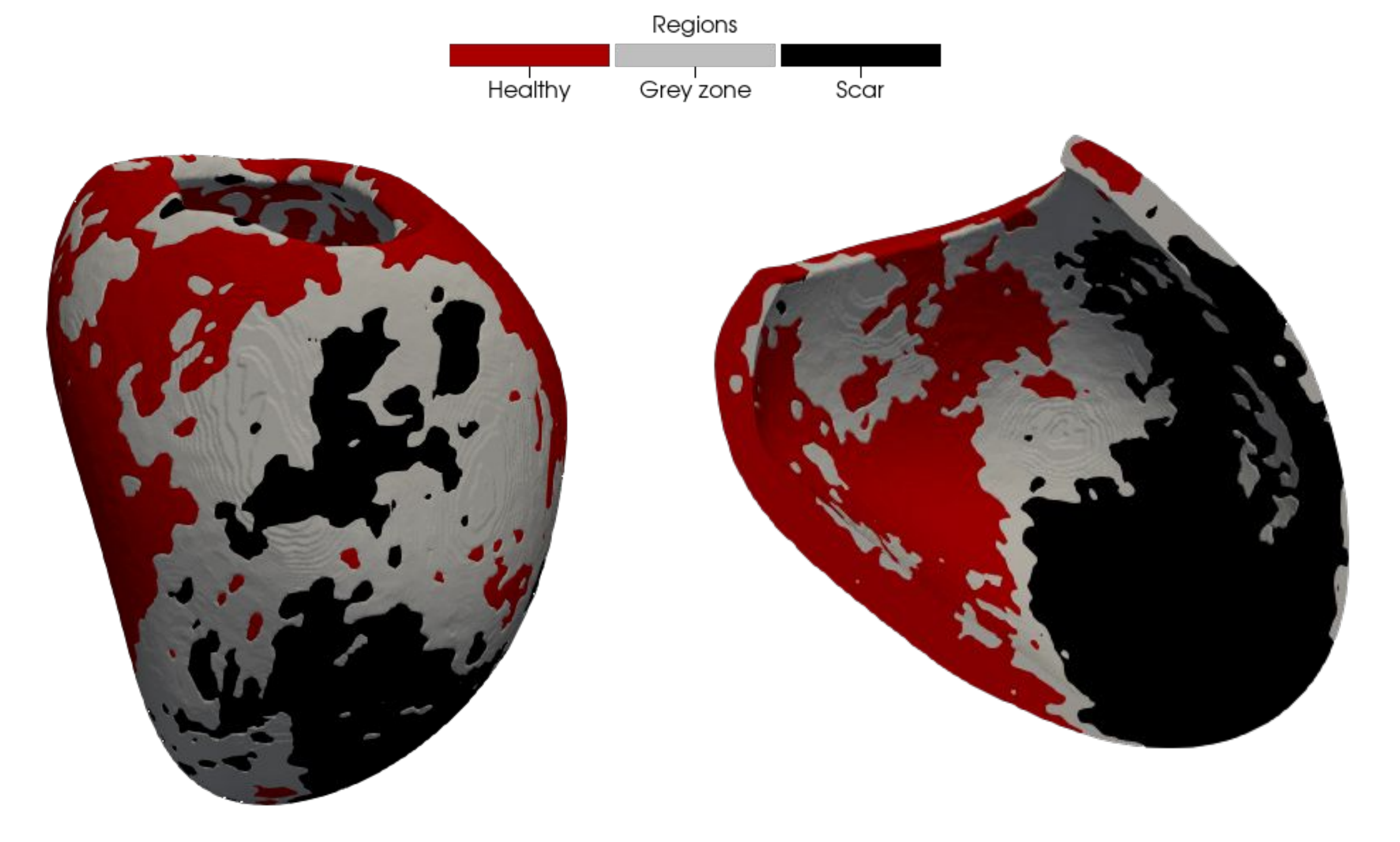}
\caption{Patient-specific LV with ICM: distribution of scars (black), grey zones (grey) and non-remodeled regions (red) over the myocardium. Volumetric view (left) and cut view (right).}
\label{fig: ischemicregion}
\end{figure}

We use a partitioned scheme to solve separately electrophysiology, activation, mechanics and cardiovascular circulation in different blocks. This segregated approach avoids to deal with a monolithic system, which is more memory demanding and entails higher computational costs, while leading to comparably accurate results \cite{Regazzoni2020PartI, Regazzoni2020PartII, Salvador2020}. For what concerns the time discretization, we use the backward differentiation formula (BDF) scheme for electrophysiology, activation and mechanics, and an explicit Runge-Kutta scheme for the circulation model.
These choices allow to accurately catch the fast time dynamics of the model variables without unbearable restrictions on the timestep, while not introducing numerical instabilities. We treat the nonlinear terms coming from both electrophysiology and activation models in a semi-implicit fashion \cite{Salvador2020}. Cardiac mechanics, which is highly nonlinear, is numerically advanced in time with a fully implicit scheme \cite{Gerbi, Salvador2020}.

We consider a first-order time splitting scheme where different timesteps are employed according to the time scale of the specific problem \cite{GerbiSegregated, GerbiPhD}. We use a fine time resolution ($\Delta t=\SI{5e-5}{\second}$) for electrophysiology, activation and circulation, whereas we use a larger timestep ($20 \Delta t$) for mechanics.

We employ the Finite Element Method (FEM) for the space discretization of electrophysiology, activation and mechanics models \cite{Quarteroni3, Quarteroni1}. We use hexahedral meshes and the $\mathbb{Q}_1$ finite element space for all core models.
This multiphysics problem presents different space resolutions according to the specific model at hand \cite{Regazzoni2020PartI, Regazzoni2020PartII}. In the framework of cardiac electrophysiology, a very fine geometrical description is required to accurately capture the electric propagation due to fine-scale phenomena arising from the continuum modeling of the cellular level, especially with the aim of reproducing and properly address arrhythmias \cite{Salvador2020}. On the other hand, cardiac mechanics requires a lower space resolution, whereas its numerical solution is more computationally demanding given the intrinsic high degree of nonlinearity: hence, using a smaller number of nodes would lead to a significant gain on the final computational time, and in particular on the assembling phase \cite{Salvador2020}.
Despite an intergrid approach, as proposed in \cite{Regazzoni2020PartI, Regazzoni2020PartII, Salvador2020}, would be beneficial, in this work we do not consider different mesh resolutions for electrophysiology and mechanics but rather use a mesh that balances numerical accuracy and computational efficiency for both models. Indeed we believe that having a consistent geometrical representation of the complex patient-specific ischemic regions, as depicted in Fig.~\ref{fig: ischemicregion}, is essential to our purposes, whereas using grids of different resolutions would introduce a mismatch among the different physical models.

For more information about the numerical details presented in this section we refer to \cite{GerbiPhD, Regazzoni2020PartI, Regazzoni2020PartII, Salvador2020}.

\section{Results}
\label{sec: numericalresults}

\begin{table}[t!]
 \center
 \caption{Patient-specific data for a LV with ICM.}
 \label{tab: clinicaldata}
 \begin{tabular}{ll}
 \hline\noalign{\smallskip}
 Parameter & Value \\
 \noalign{\smallskip}\hline\noalign{\smallskip}
 EF & 13 $\%$ \\
 SV & 39 mL \\
 HR & 57 \\
 MAP & 80 mmHg \\
 SBP & 110 mmHg \\
 DBP & 69 mmHg \\
 \noalign{\smallskip}\hline
 \end{tabular}
\end{table}

\begin{figure}[t!]
\centering
\includegraphics[keepaspectratio, width=1.0\textwidth]{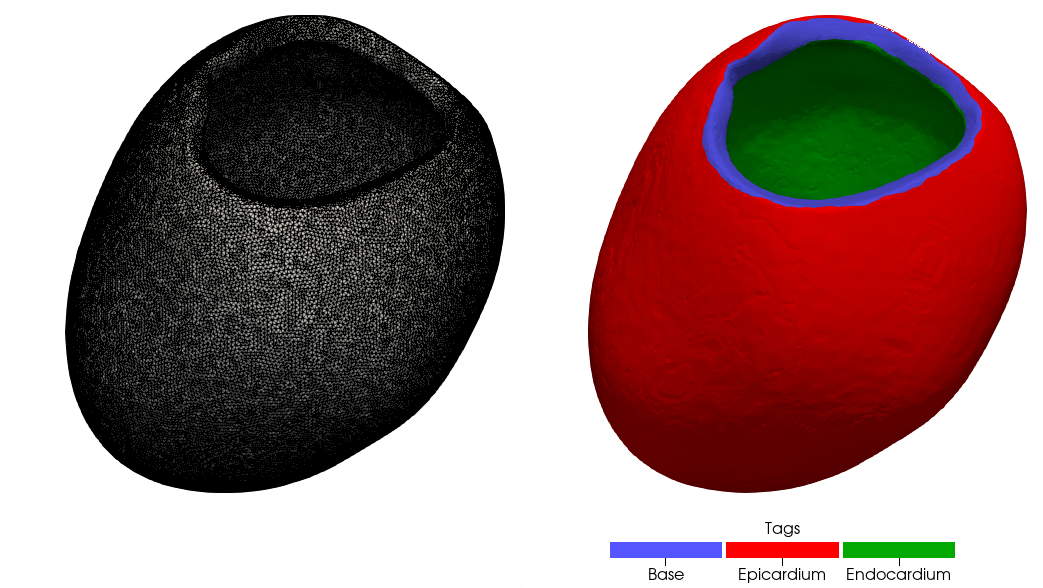}
\caption{Initial tetrahedral mesh (left) and distribution of tags over the myocardium (right).}
\label{fig: mesh_tags}
\end{figure}

\begin{figure}[t!]
\centering
\includegraphics[keepaspectratio, width=1.0\textwidth]{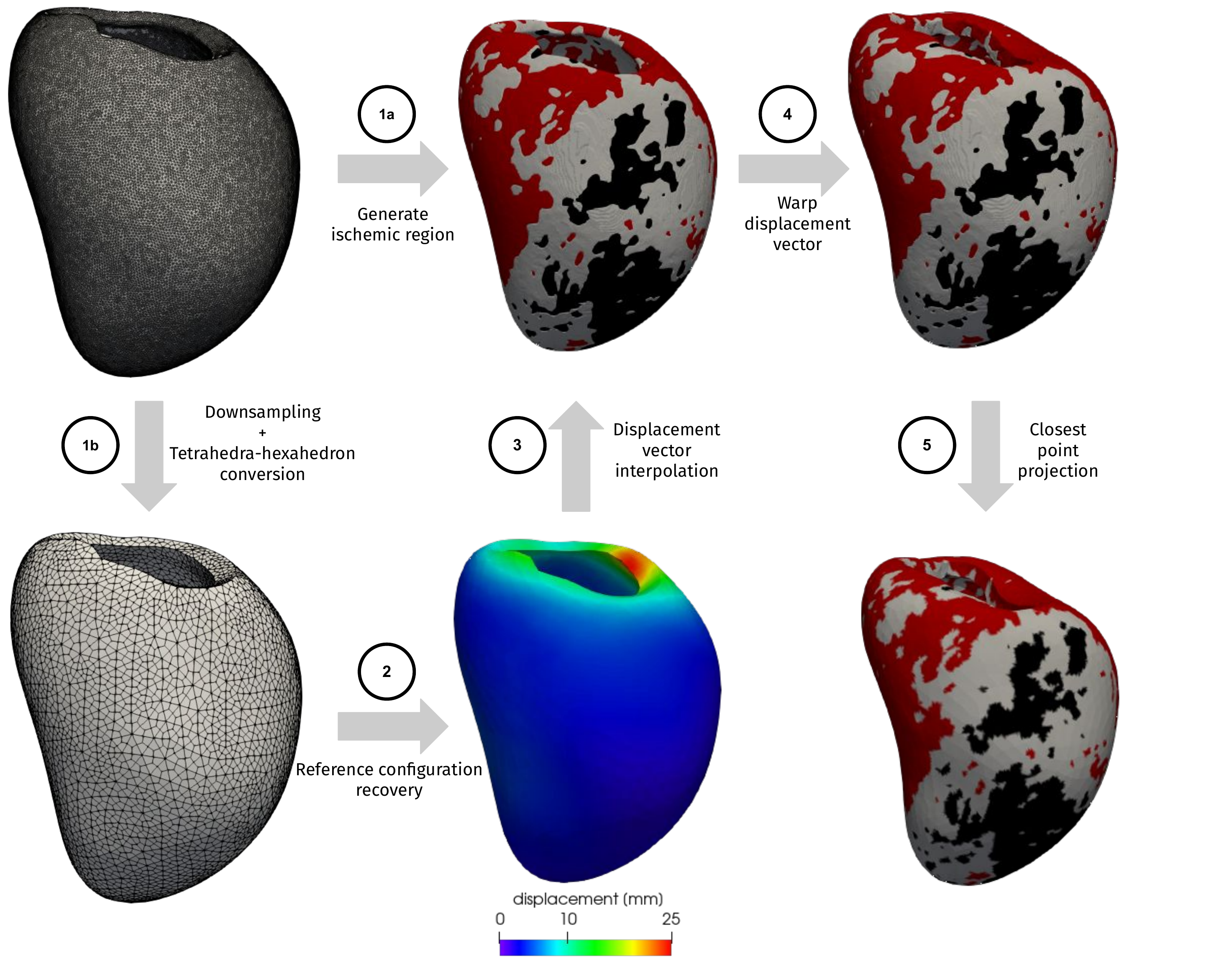}
\caption{Preprocessing pipeline for a patient-specific LV with ICM. The ischemic regions are mapped on the high-resolution tetrahedral mesh (step 1a). Then, downsampling and conversion to hexahedral elements are performed (step 1b). We compute the reference configuration on the latter mesh (step 2). The displacement vector is accurately and efficiently interpolated on the high-resolution mesh with the ischemic regions, by using the interpolant proposed in \cite{Regazzoni2020PartII} (step 3). Finally, we warp the tetrahedral mesh with the interpolated reference configuration displacement (step 4) and we perform closest point projection of the ischemic regions distribution on the hexahedral mesh used in our numerical simulations (step 5).}
\label{fig: ischemic_region_interpolation}
\end{figure}

\begin{figure}[t!]
\centering
\includegraphics[keepaspectratio, width=1.0\textwidth]{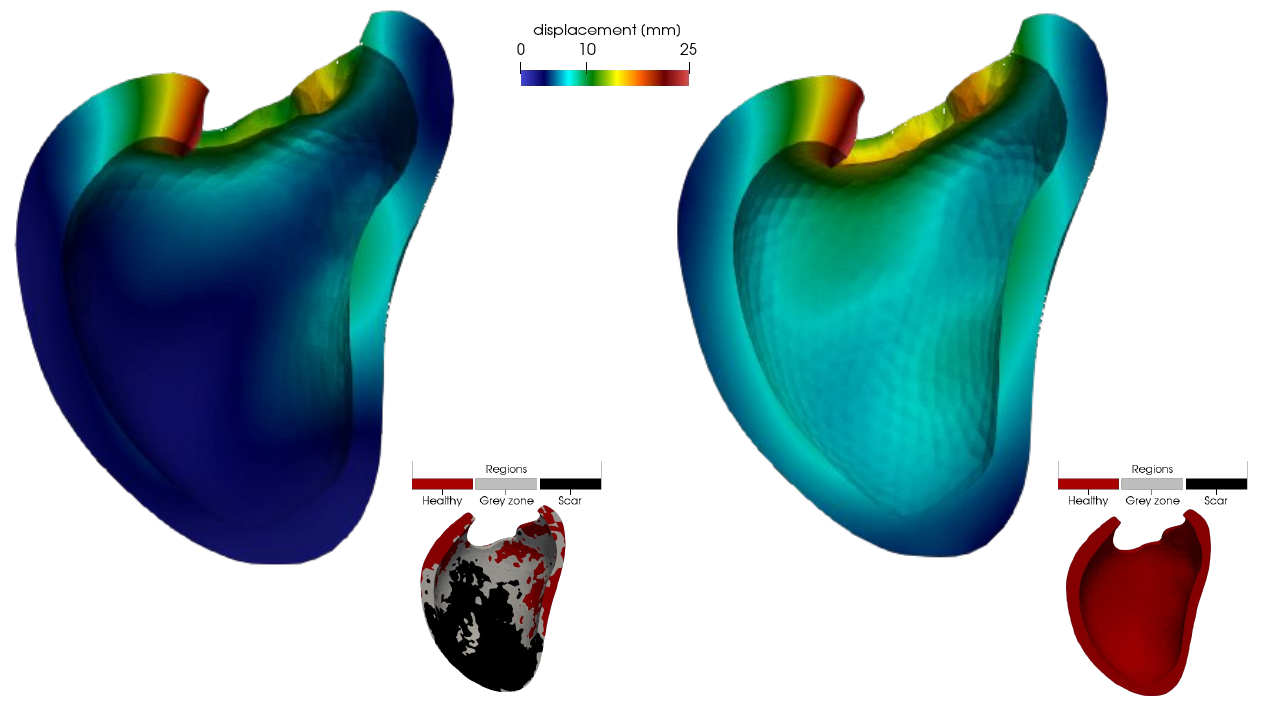}
\caption{Reference configuration recovery for a patient-specific LV with ICM. We depict the displacement to reach the unstressed geometry $\Omega_0$ in case the ischemic regions distribution is imported (left) or when fully healthy conditions are considered (right). Scars, which are mostly localized at the apex, prevent the unloading for part of the myocardium.}
\label{fig: referenceconfiguration}
\end{figure}

We depict in Fig.~\ref{fig: ischemicregion} the geometric model of the patient-specific LV. The distribution of the ischemic regions is mapped from the LGE-MRI. We use the Segment cardiac image analysis software package \cite{Segment} to segment the Cine MRI of the patient. From this segmentation, we retrieve the volume of the LV blood pool over time, the stroke volume (SV), the ejection fraction (EF) and the heart rate (HR). We also collected some pressure data, namely the mean arterial pressure (MAP), the systolic blood pressure (SBP) and the diastolic blood pressure (DBP). All these clinical data are shown in Tab. \ref{tab: clinicaldata}.

The geometry is tetrahedralized using elements of average edge length equal to 0.35 mm. In Fig.~\ref{fig: mesh_tags} we depict the tetrahedral mesh and the tags for base, epicardium and endocardium. Meshing and ischemic region mapping from imaging are performed with the commercial software Materialise Mimics \cite{Mimics}.
The VMTK software -- the Vascular Modelling Toolkit \cite{VMTK, Antiga, Fedele2021} -- is used to generate the corresponding hexahedral mesh and to downsample it to the desired resolution employed in our numerical simulations (average element diameter $h_\text{mean}$=1.5 mm, 638'048 elements and 692'535 vertices). We also use VMTK to accurately map the ischemic region distribution from the original high resolution mesh to the downsampled one, after computing the reference configuration, as shown in Fig.~\ref{fig: ischemic_region_interpolation}, by using the technique proposed in \cite{Regazzoni2020PartII}.

We use the Bayer-Blake-Plank-Trayanova algorithm \cite{Bayer, Piersanti2020} to generate the fibers distribution (field $\fZero$) for our geometry, with $\alpha_{\text{epi}} = -60^{\circ}$, $\alpha_{\text{endo}} = 60^{\circ}$, $\beta_{\text{epi}} = 20^{\circ}$ and $\beta_{\text{endo}} = -20^{\circ}$.

We compute the reference configuration of this patient-specific LV by solving the inverse problem proposed in \cite{Regazzoni2020PartI, Regazzoni2020PartII}. In this way, we get the geometry in a stress-free condition, as if the blood was completely removed from the LV. Differently from \cite{Regazzoni2020PartI, Regazzoni2020PartII}, the new formulation (Eqs.~\eqref{eqn: monodomainionic},~\eqref{eqn: activestrain},~\eqref{eqn: mechanics},~\eqref{eqn: Guccione_pt1},~\eqref{eqn: Guccione_pt2}) accounting for the heterogeneous coefficient $\eta = \eta(\boldsymbol{x})$, has a clear influence on the numerical solution to this inverse problem, as shown in Fig.~\ref{fig: referenceconfiguration}.

We present numerical results in the context of cardiac electromechanics, both in SR and during VT.
First, we perform parameters calibration in SR to fit the available clinical data. For this stage, we employ our 3D electromechanics model and a simplified 0D afterload model \cite{GerbiPhD, Salvador2020}.
Then, for the VT, we consider the coupling between the 3D electromechanics model and the 0D circulation model. A calibration of the fully coupled model is out scope and would be very challenging, especially due to the lack of data for this patient on atria, right ventricle and pulmonary circulation. Moreover, our purpose is to show the capabilities of this framework to investigate VTs in complex patient-specific ventricles, by providing both quantitative and qualitative observations on arrhythmias when electrophysiology, mechanics and fluid dynamics are coupled together. Further details about the latter aspect are provided in Sec. \ref{sec: arrhythmiasimulations}.

The mathematical models of Sec.~\ref{sec: mathematicalmodeling} and the numerical methods of Sec.~\ref{sec: numericaldiscretization} have been implemented in \texttt{life\textsuperscript{x}} (\url{https://lifex.gitlab.io/lifex}), a high-performance \texttt{C++} library developed within the iHEART project and based on the \texttt{deal.II} (\url{https://www.dealii.org}) Finite Element core \cite{dealII92}.

The numerical simulations were performed on a HPC facility available at MOX for the iHEART project. The entire cluster is endowed with 8 Intel Xeon Platinum 8160 processors, for a total of 192 computational cores and a total amount of 1.5TB of available RAM.

\subsection{Simulations in sinus rhythm}
\label{sec: sinusrhythmsimulations}

\newcommand{\EMsnapshotSR}[2]{
	\subfloat[][$t = \SI{#2}{\second}$]{\includegraphics[width=0.6\textwidth]{pictures/EM_simulation_SR/EM_simulation_SR_#1.png}}}

\begin{figure}[t!]
	\advance\leftskip-1.4cm
	\captionsetup[subfigure]{labelformat=empty}
	\EMsnapshotSR{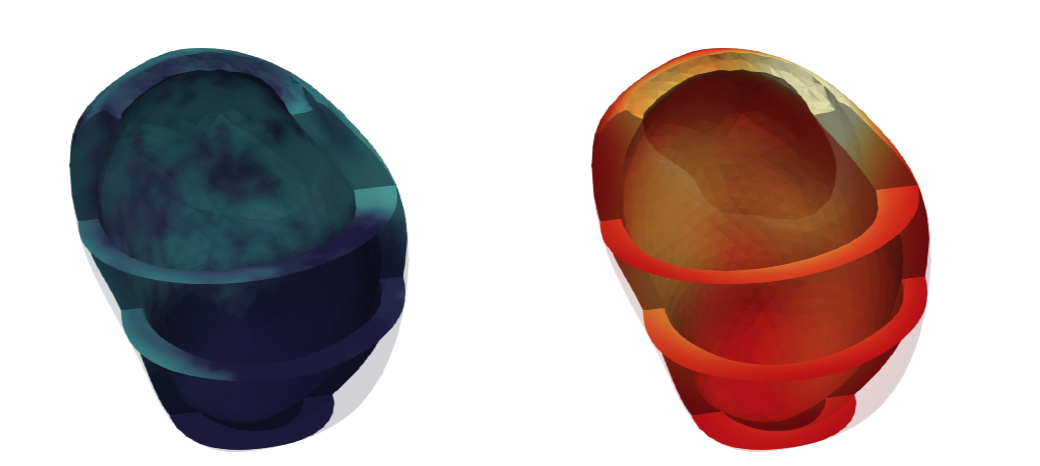}{0.03} \quad
	\EMsnapshotSR{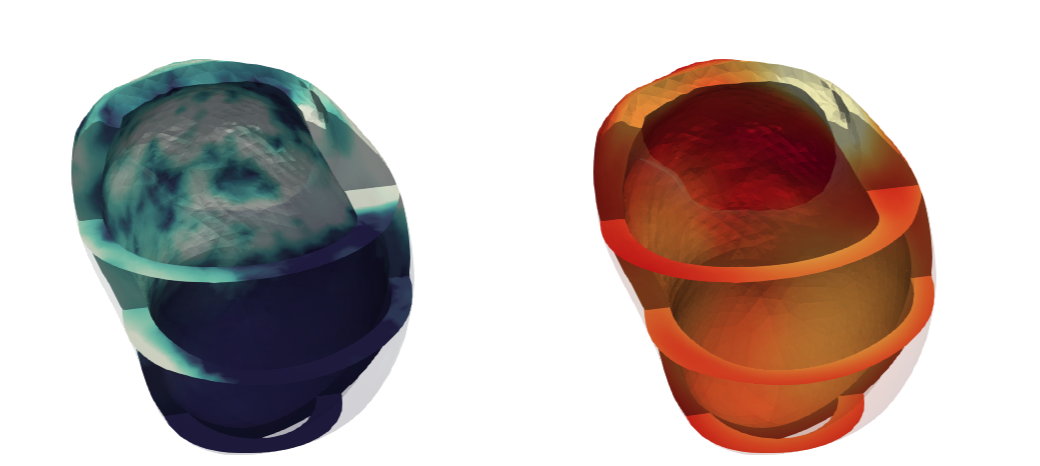}{0.06} \\
	\EMsnapshotSR{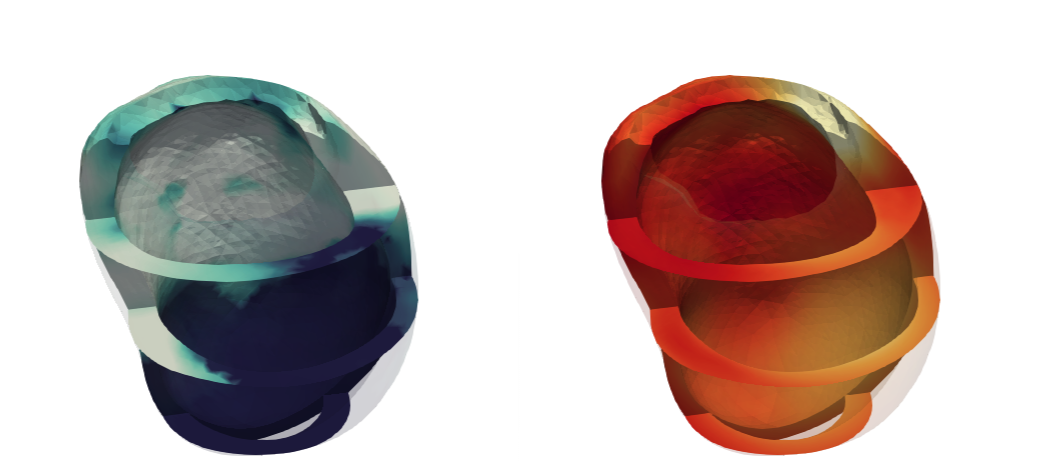}{0.20} \quad
	\EMsnapshotSR{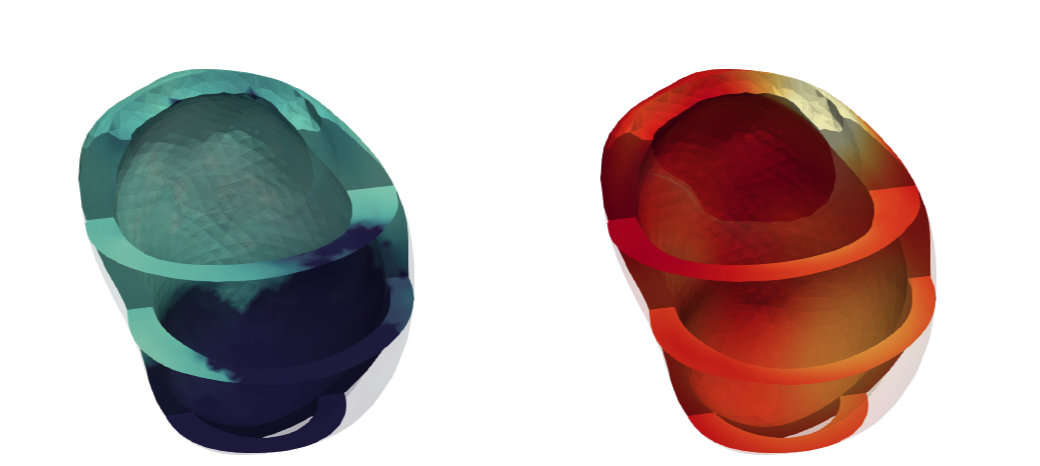}{0.25} \\
	\EMsnapshotSR{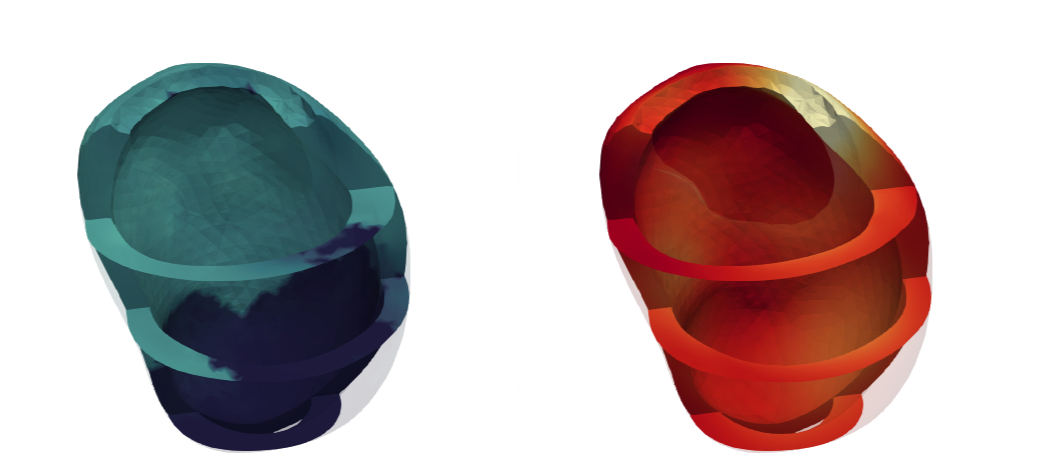}{0.30} \quad
	\EMsnapshotSR{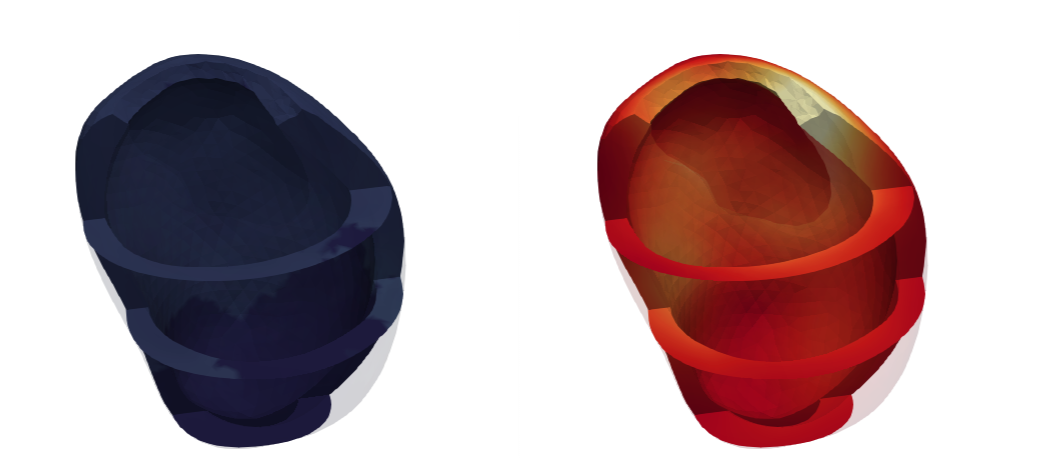}{0.40} \\
	\EMsnapshotSR{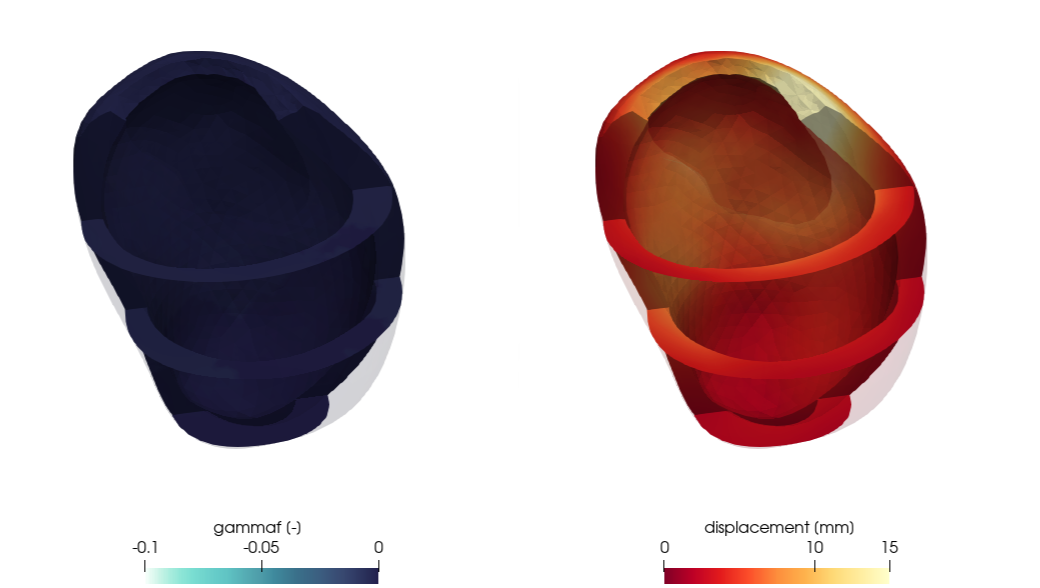}{0.60} \quad
	\EMsnapshotSR{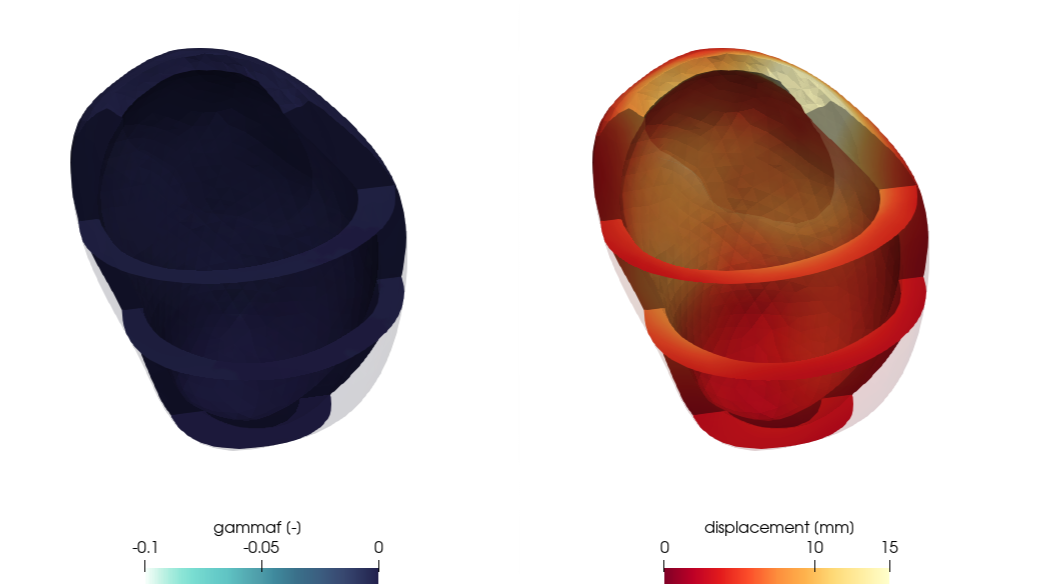}{0.90}

	\caption{Evolution in SR of $\gamma_f$ and displacement magnitude $|\displ|$ for a patient-specific LV with ICM. Each picture is warped by the displacement vector.}
	\label{fig: EM_simulation_SR}
\end{figure}

\begin{figure}[t!]
\centering
\advance\leftskip-1.2cm
\includegraphics[keepaspectratio, width=1.2\textwidth]{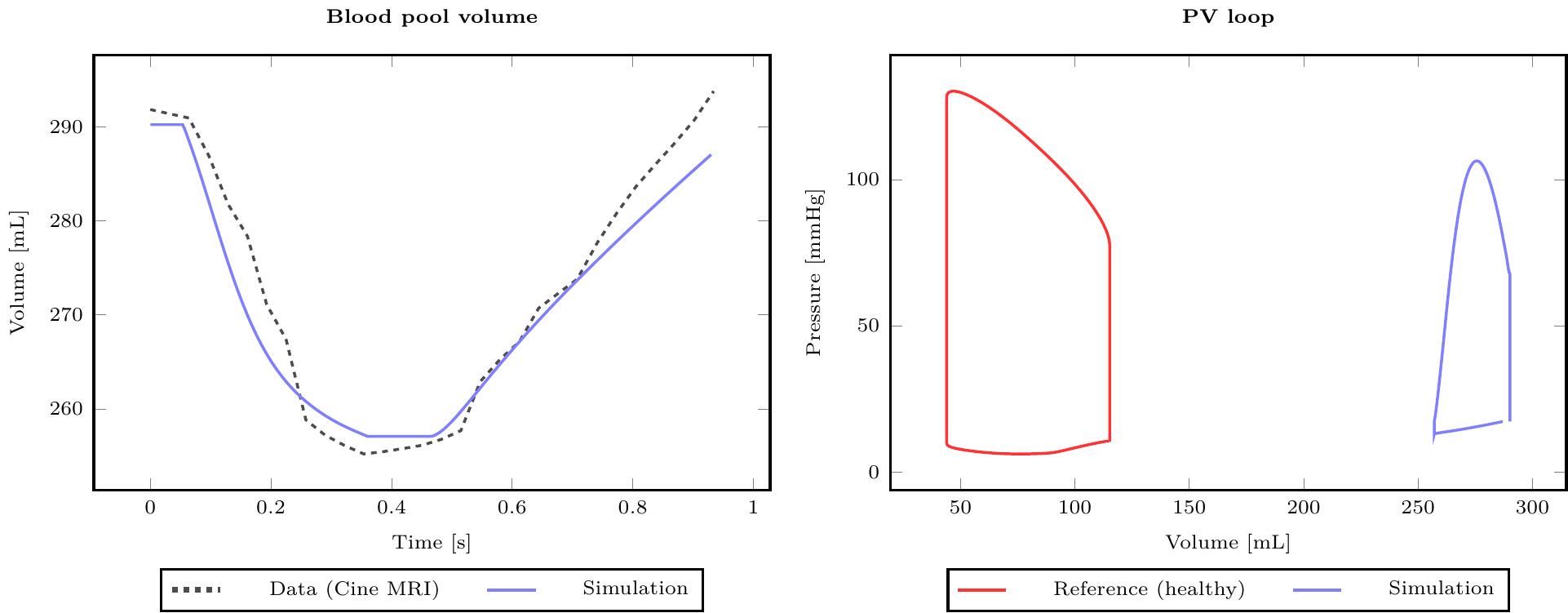}
\caption{Electromechanics simulation of a patient-specific LV with ICM: blood pool volume over time (left) and PV loop (right). We highlight the comparison with clinical data and a reference healthy LV (taken from \cite{Regazzoni2020PartII}).}
\label{fig: pt339_volume_pvloop}
\end{figure}

We run several numerical simulations in SR to perform the calibration of the parameters of our electromechanical model. In particular, we tuned $\muA$ for active strain, $\BCmecKepiN$, $\BCmecKepiT$, $\BCmecCepiN$ and $\BCmecCepiT$ for mechanics, resistance $R$ and capacitance $C$ of the 2-element Windkessel model.
We do not calibrate the passive mechanics parameters related to the Guccione constitutive law. Indeed, we are unable to fit a Klotz curve for this patient-specific case \cite{Klotz}. We reuse the same parametrization provided in \cite{Sack} on failing swine hearts.
Moreover, even if we do not perform a quantitative strain analysis, we get a good qualitative match between the displacement field observed from our numerical simulation and the one coming from Cine MRI.
We remind that the reference configuration recovery is repeated for each new set of parameters.
We apply a current $\EPIappReduced(\boldsymbol{x}, t)$, with a cubic distribution in space and peak $\EPIappReducedMax$, for a duration of $\EPIappDuration$, in different regions of the myocardium, to trigger the electrical signal in the LV. While the Purkinje network is not explicitly represented here \cite{Vergara, Romero}, the stimuli generate a realistic propagation pattern.
We report in Appendix~\ref{app:params} the final configuration that we used for the SR simulation. In Fig.~\ref{fig: EM_simulation_SR} we depict the evolution of $\gamma_f$ and displacement magnitude $|\displ|$ over different timesteps of the simulation. We highlight the heterogeneity in the activation of the myocardium, that leads to different contractility according to the specific region. Dense scars, which approximately occupy half of the LV, do not contract. Displacement magnitude $|\displ|$ is dominant at the base, so that the LV is pushed towards the apex.
In Fig.~\ref{fig: pt339_volume_pvloop}, we show the evolution of the LV pressure and volume for the patient. We provide a good quantitative match with Cine MRI for the evolution of blood pool volume over time, as well as for all values reported in Tab.~\ref{tab: clinicaldata}. Moreover, we underline the major differences with respect to a reference healthy PV loop (baseline simulation in \cite{Regazzoni2020PartII}). The presence of ICM significantly dilates the LV over the years. As a consequence, the PV loop moves to the right and both the end diastolic volume (EDV) and the end systolic volume (ESV) increase. The end diastolic pressure (EDP) is generally higher with ICM, whereas the pressure peak is lower. Indeed, the capability of the LV to push blood into the aorta is strongly impaired. This coincides with smaller contractility and reduced EF.

\subsection{Simulations under arrhythmia}
\label{sec: arrhythmiasimulations}

\begin{figure}[t!]
\centering
\includegraphics[keepaspectratio, width=0.6\textwidth]{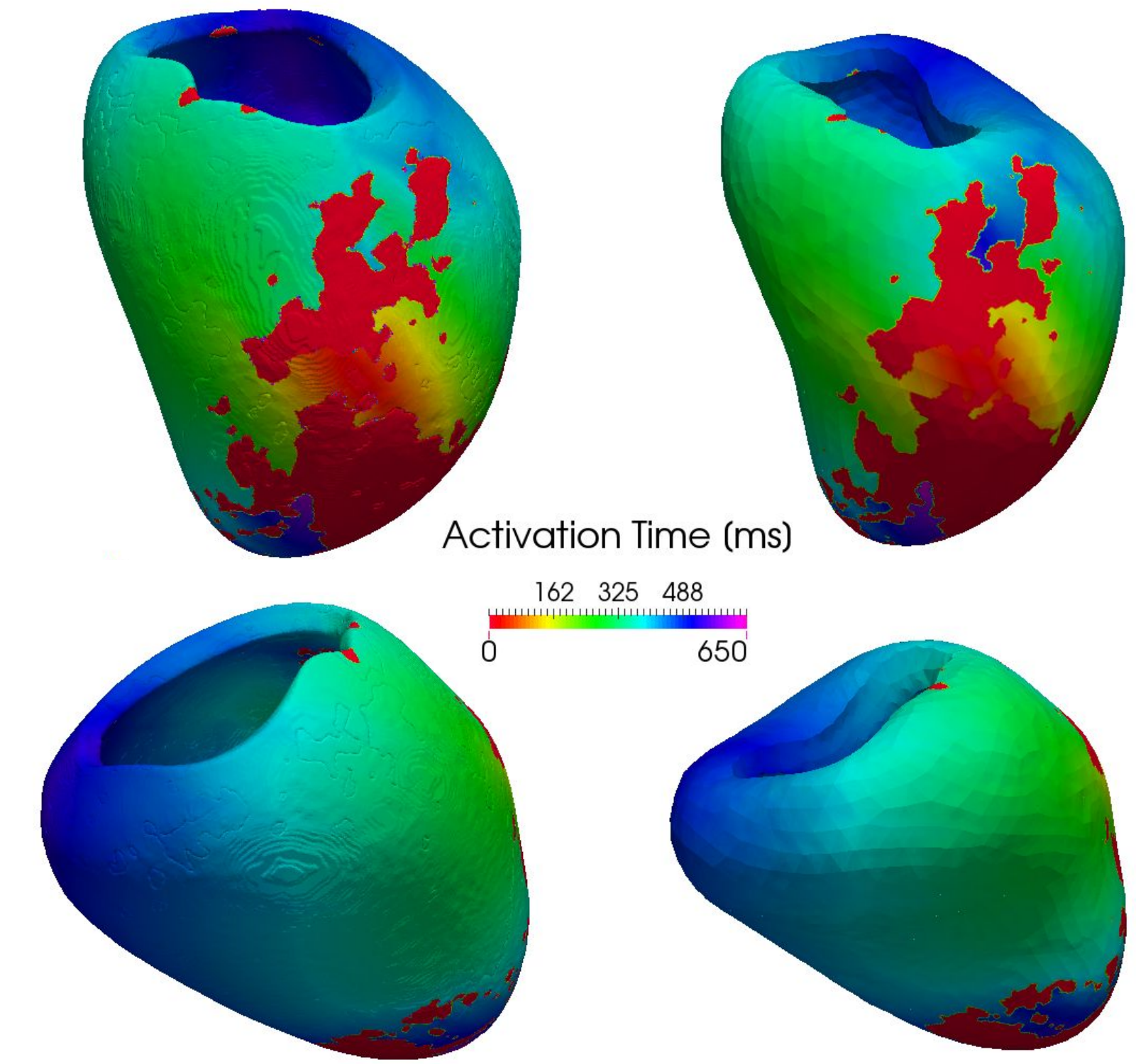}
\caption{Activation time for a patient-specific LV with ICM. Simulations from CARP (left) and \texttt{life\textsuperscript{x}} (right). We use a tetrahedral mesh ($h_\text{mean}=0.35$ mm) in CARP. We consider the reference configuration, meshed with hexahedral elements ($h_\text{mean}=1.5$ mm), in \texttt{life\textsuperscript{x}}. After tuning the conductivities of the monodomain equation properly, minor differences in terms of activation times can be observed between the two cases.}
\label{fig: pt339_AT}
\end{figure}

\newcommand{\EMsnapshotVT}[2]{
	\subfloat[][$t = \SI{#2}{\second}$]{\includegraphics[width=0.55\textwidth]{pictures/EM_simulation_VT_1/EM_simulation_VT_#1.png}}}

\begin{figure}[p]
	\advance\leftskip-0.9cm
	\captionsetup[subfigure]{labelformat=empty}
	\EMsnapshotVT{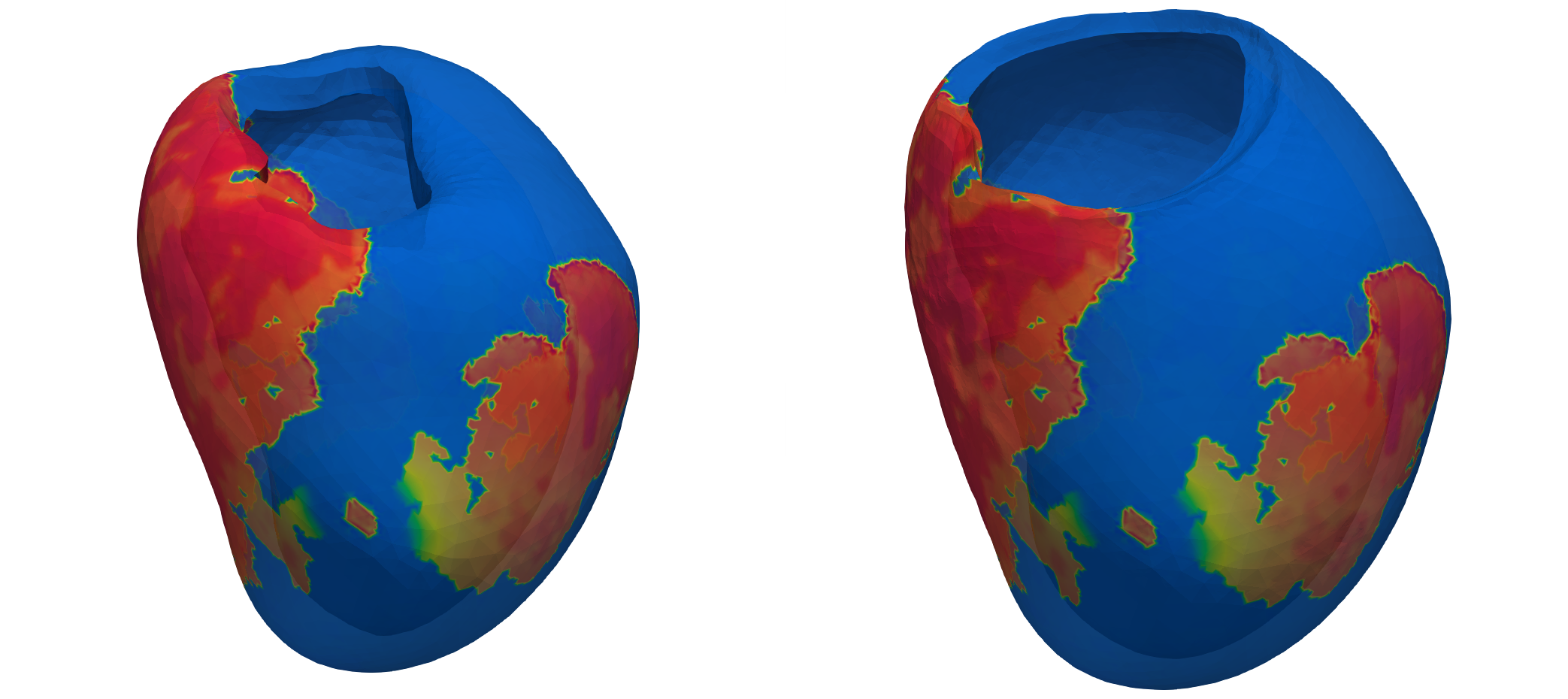}{0.5} \quad
	\EMsnapshotVT{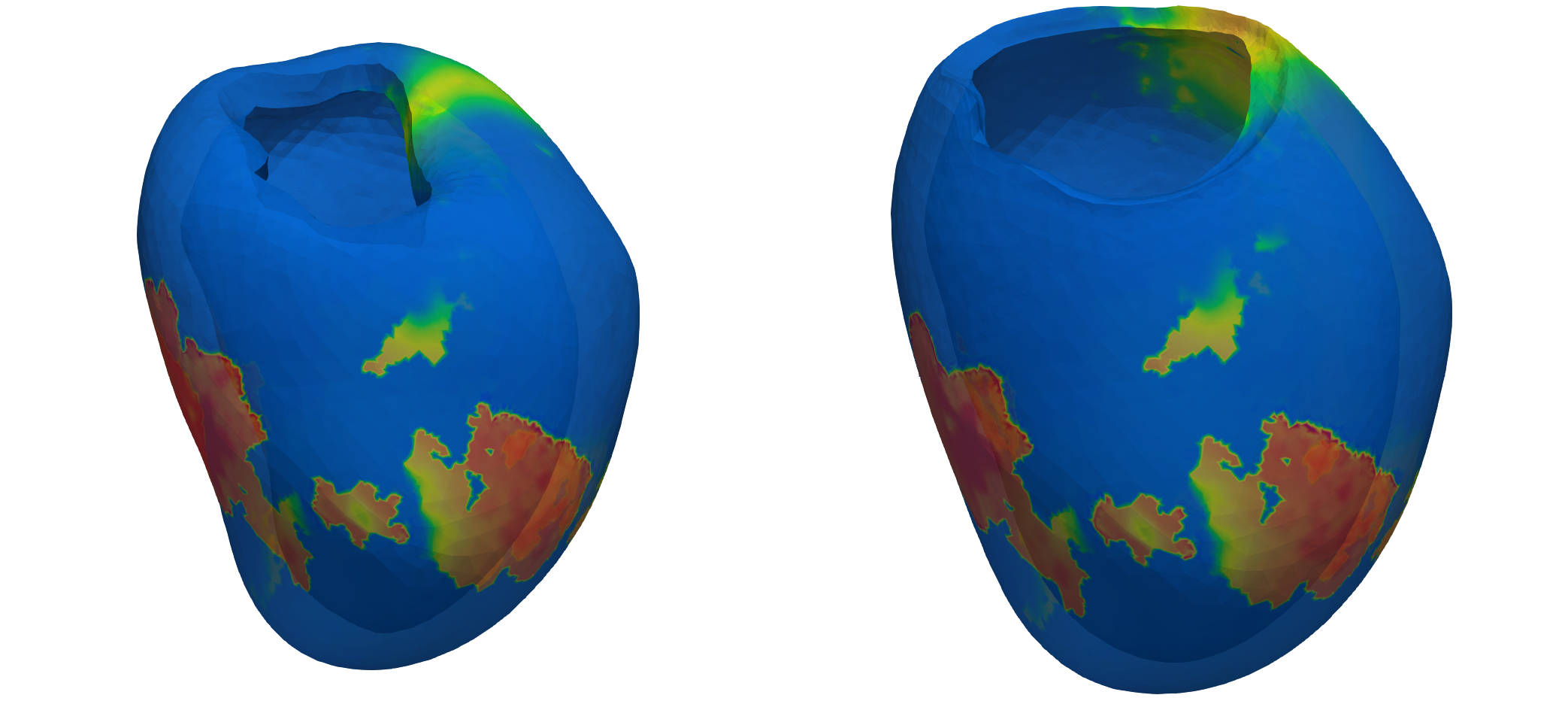}{1.0} \\
	\EMsnapshotVT{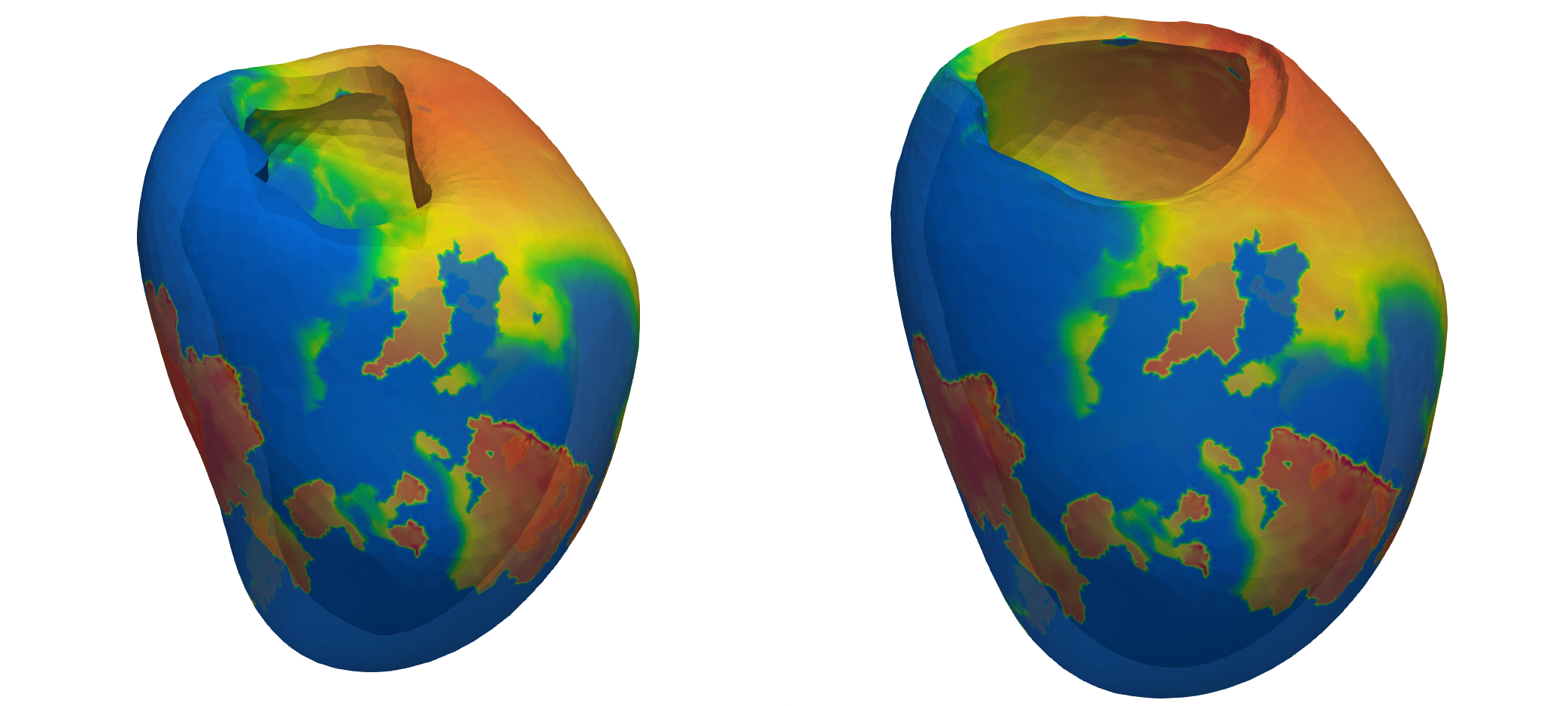}{1.5} \quad
	\EMsnapshotVT{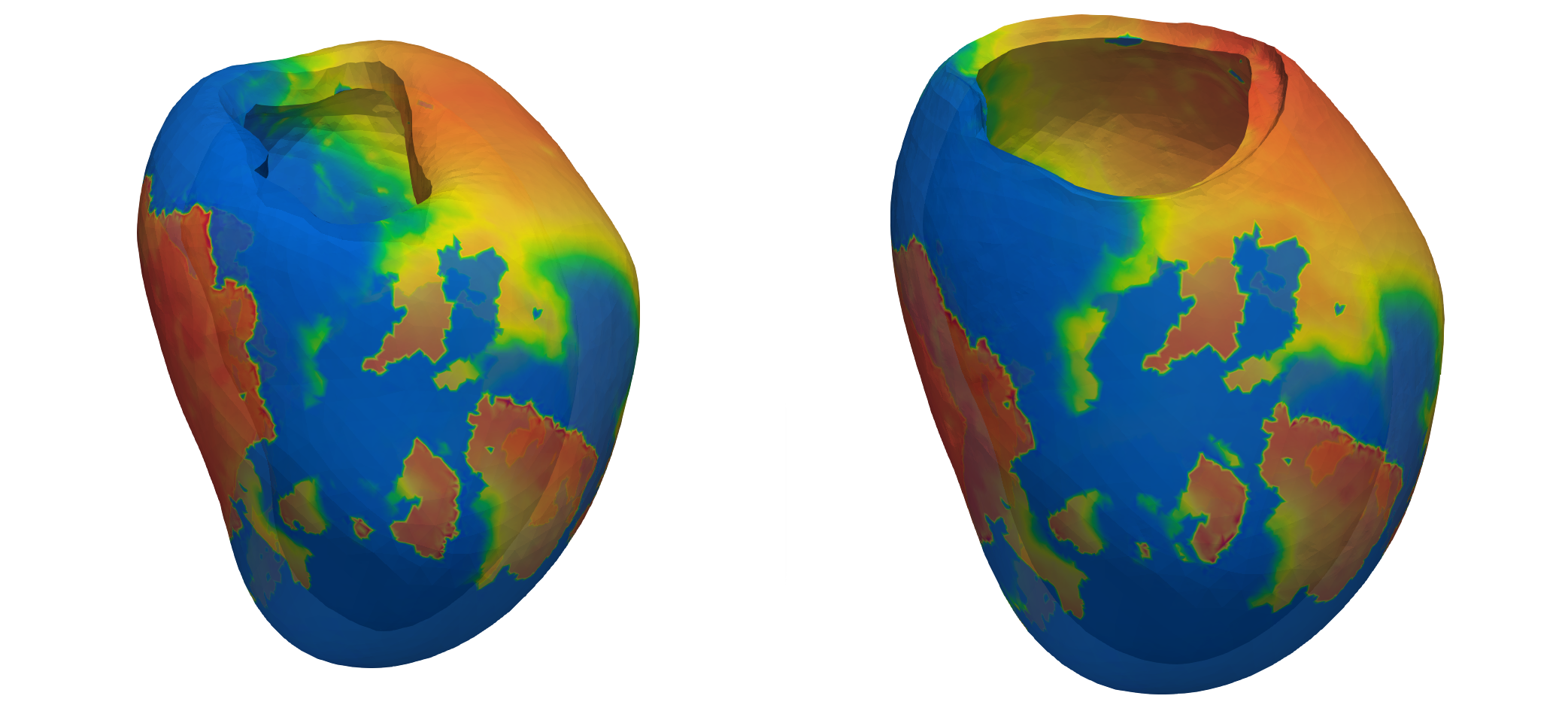}{2.0} \\
	\EMsnapshotVT{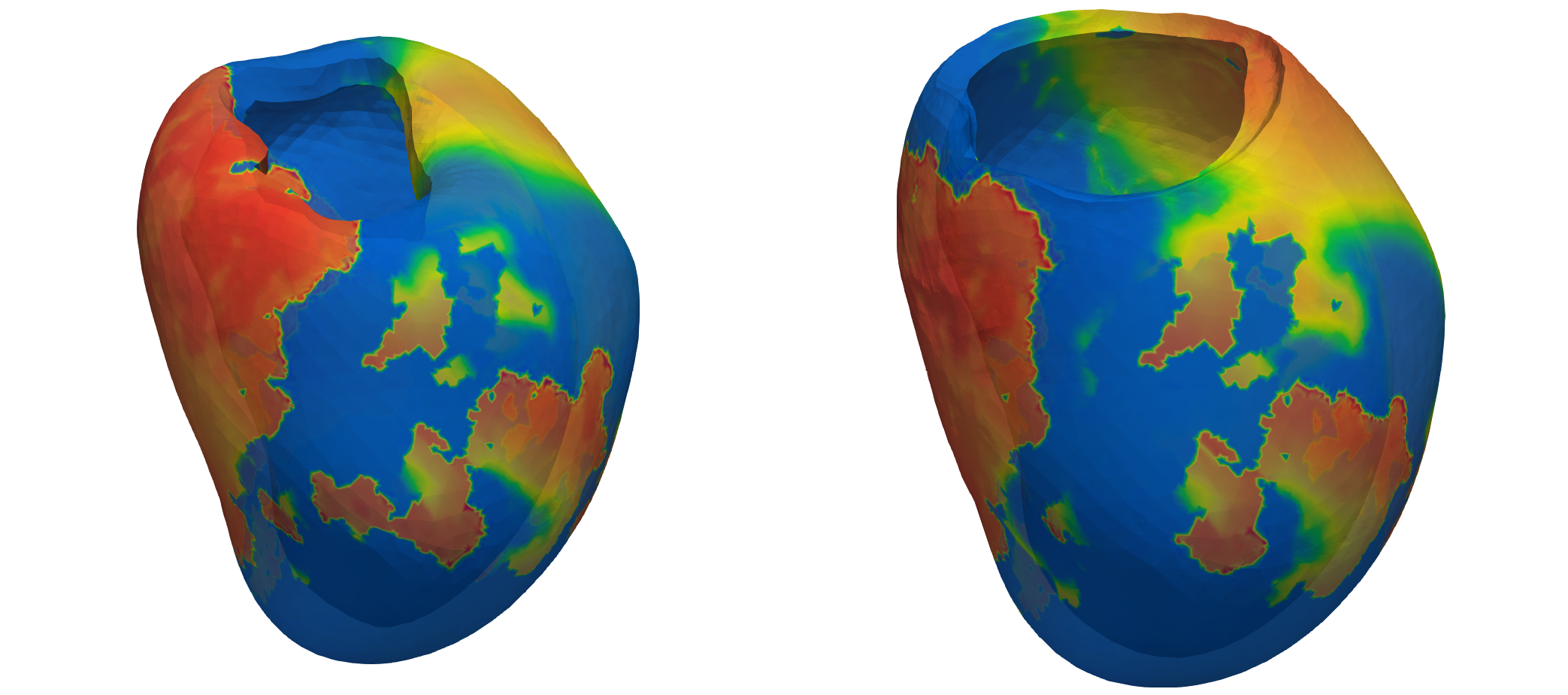}{2.5} \quad
	\EMsnapshotVT{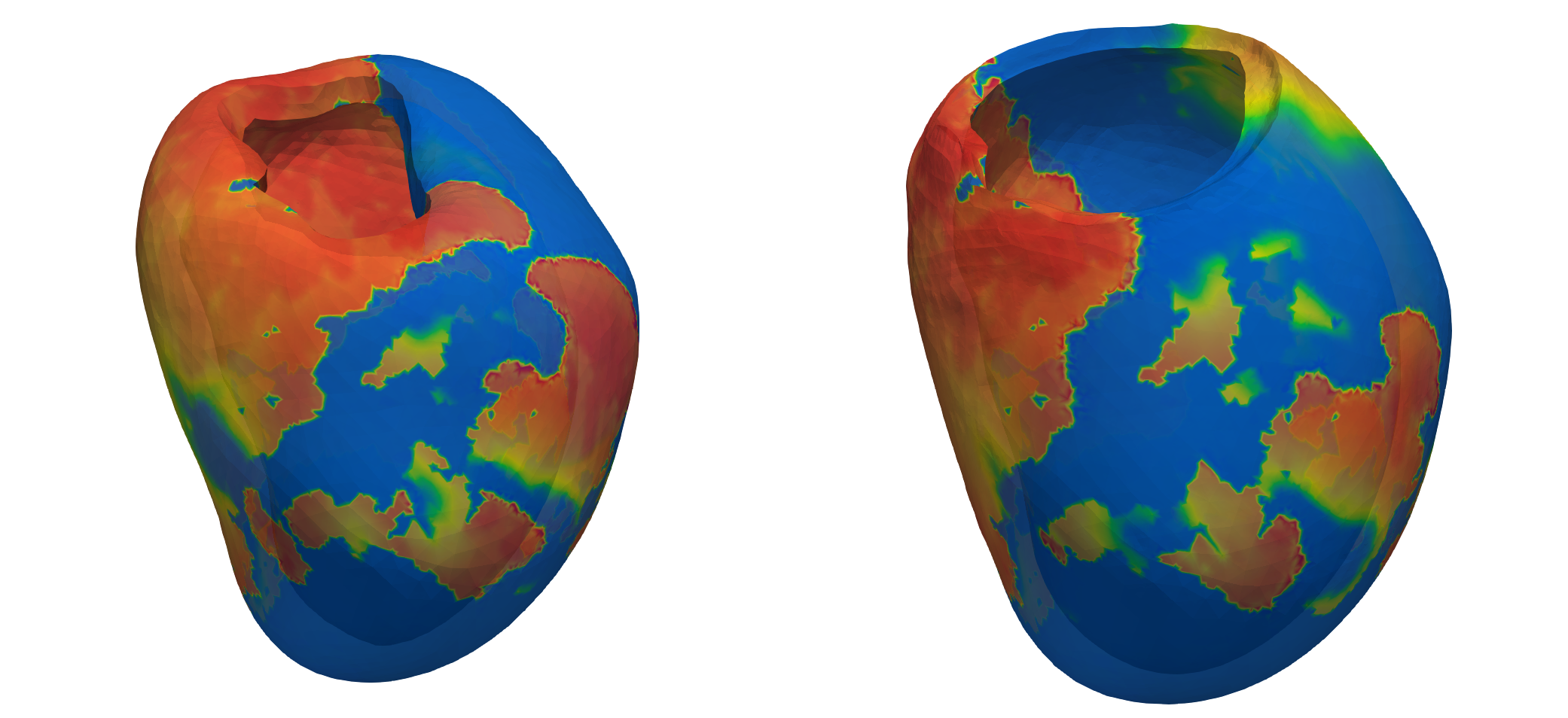}{3.0} \\
	\EMsnapshotVT{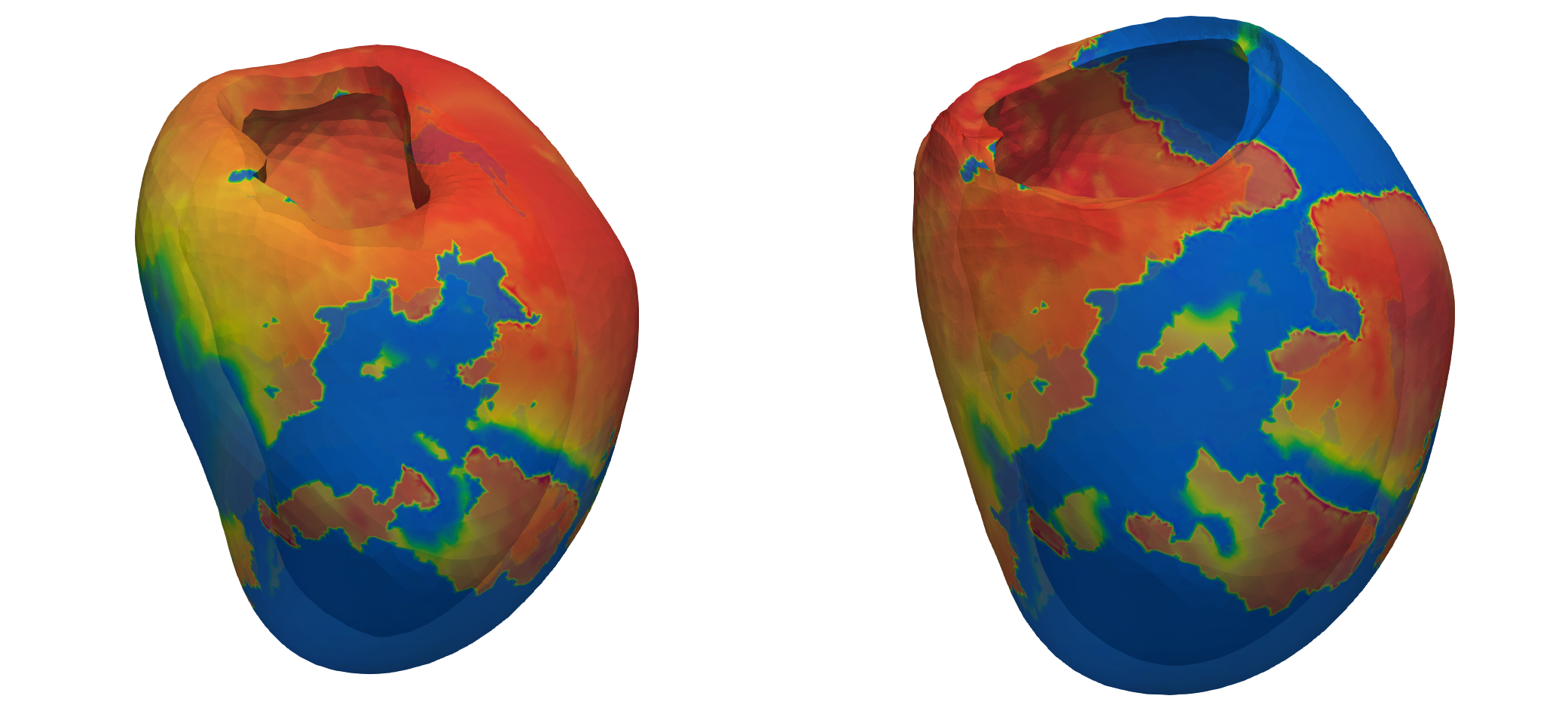}{3.5} \quad
	\EMsnapshotVT{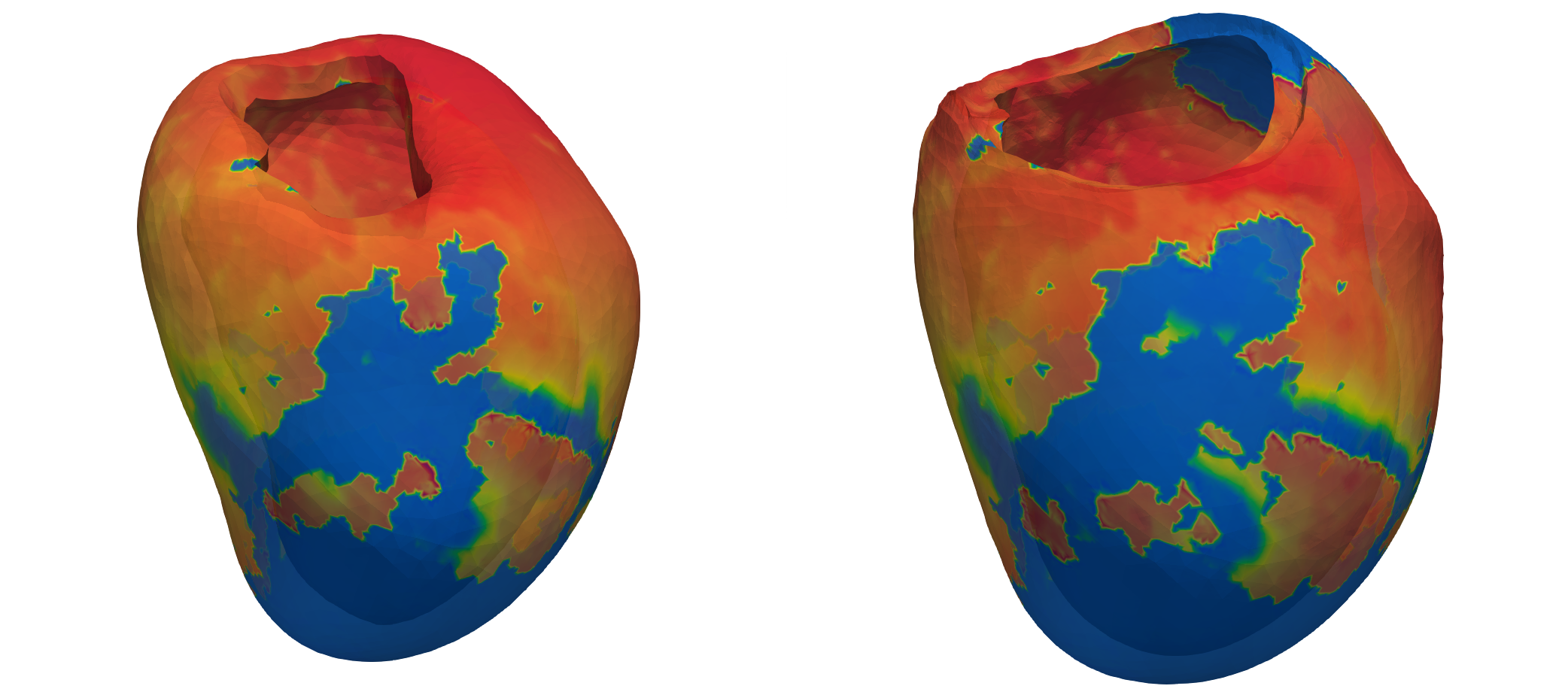}{4.0} \\
	\EMsnapshotVT{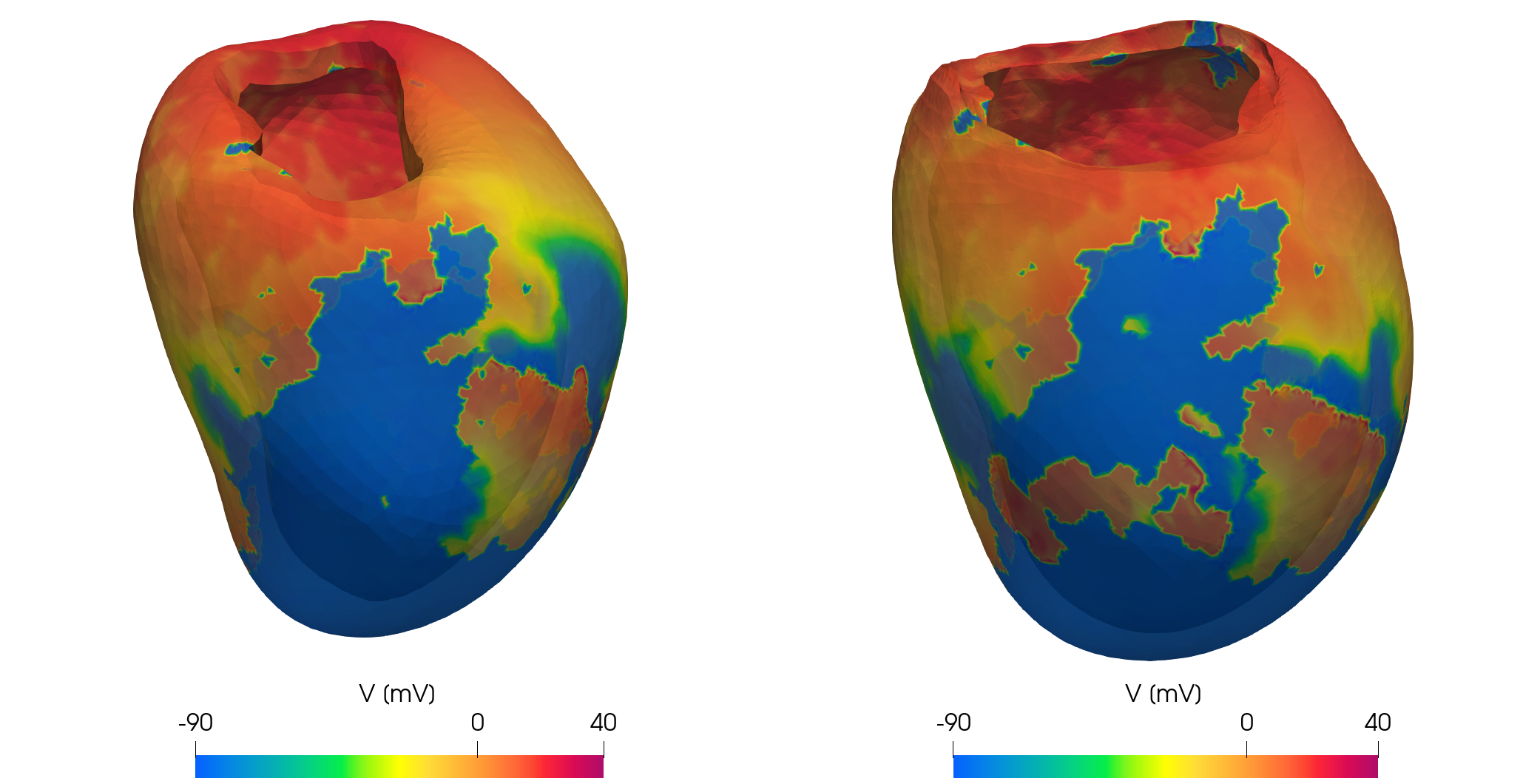}{4.5} \quad
	\EMsnapshotVT{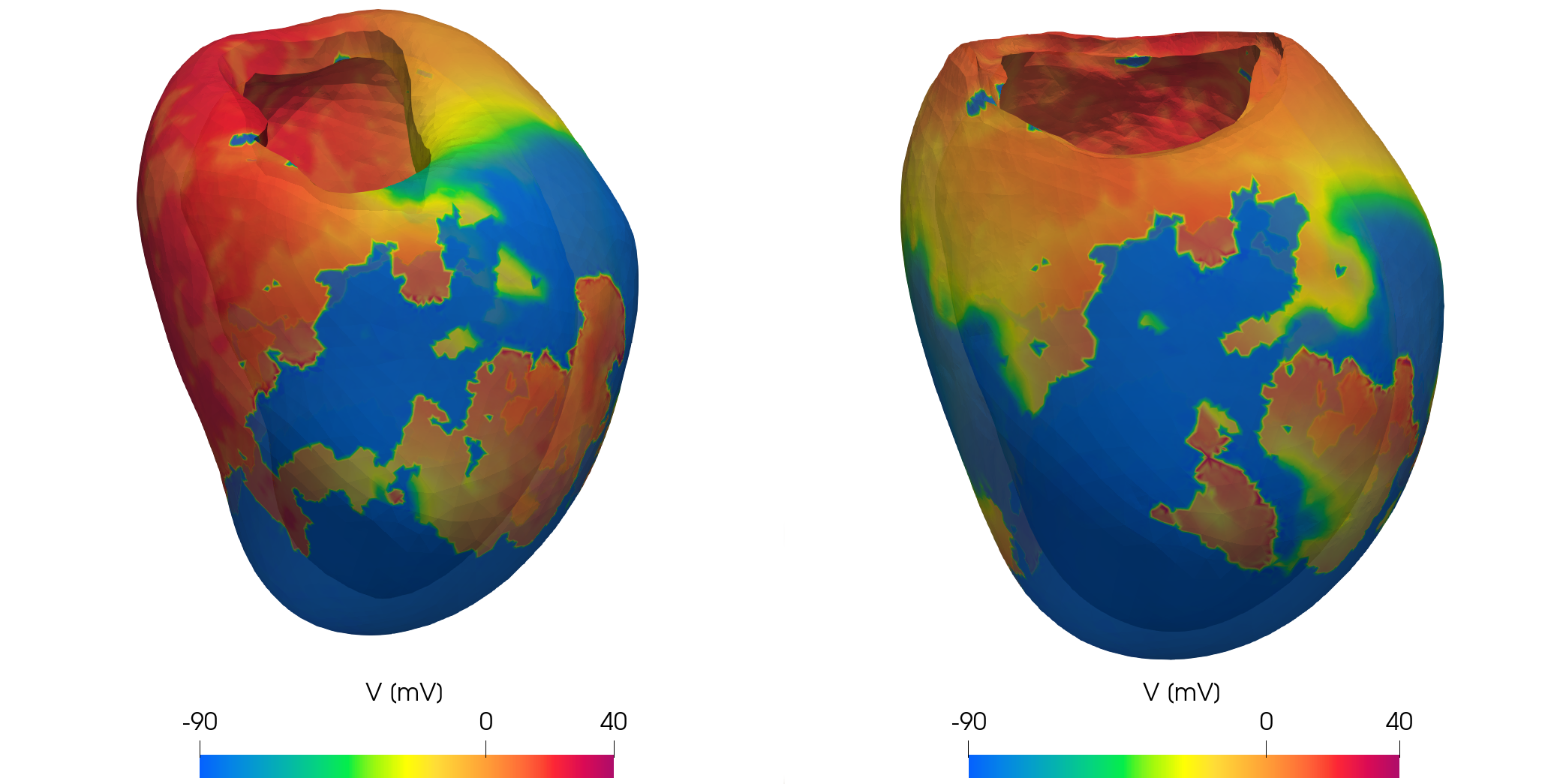}{5.0}

	\caption{Propagation of the transmembrane potential $V=85.7 u - 84$ during VT for a patient-specific LV with ICM. Electrophysiology simulation (left) runs on the geometry retrieved from LGE-MRI (\textit{i.e.} without reference configuration recovery). Electromechanics simulation (right) is warped by the displacement vector.}
	\label{fig: EM_baseline_simulation_V_displacement_1}
\end{figure}

\newcommand{\EMsnapshotVTVT}[2]{
	\subfloat[][$t = \SI{#2}{\second}$]{\includegraphics[width=0.55\textwidth]{pictures/EM_simulation_VT_2/EM_simulation_VT_#1.png}}}

\begin{figure}[p]
	\advance\leftskip-0.9cm
	\captionsetup[subfigure]{labelformat=empty}
	\EMsnapshotVTVT{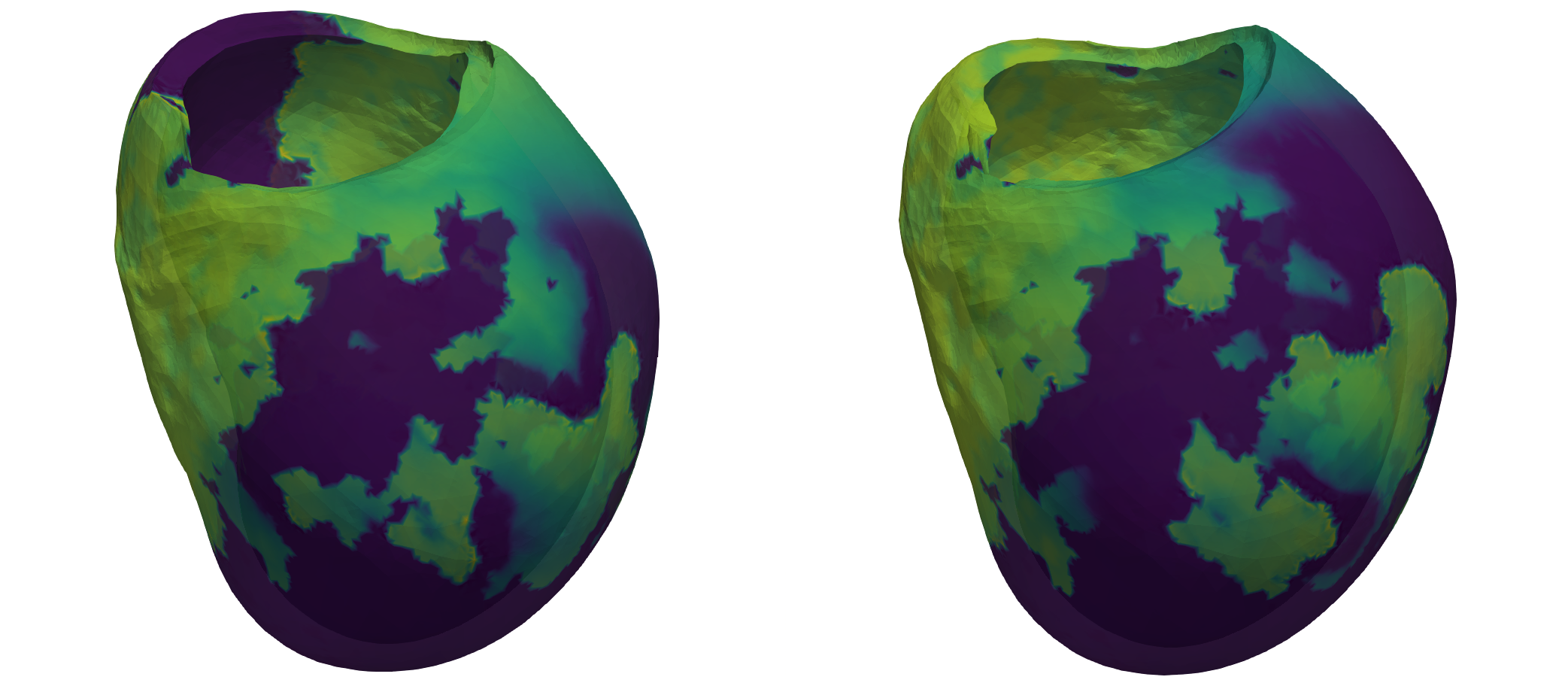}{5.5} \quad
	\EMsnapshotVTVT{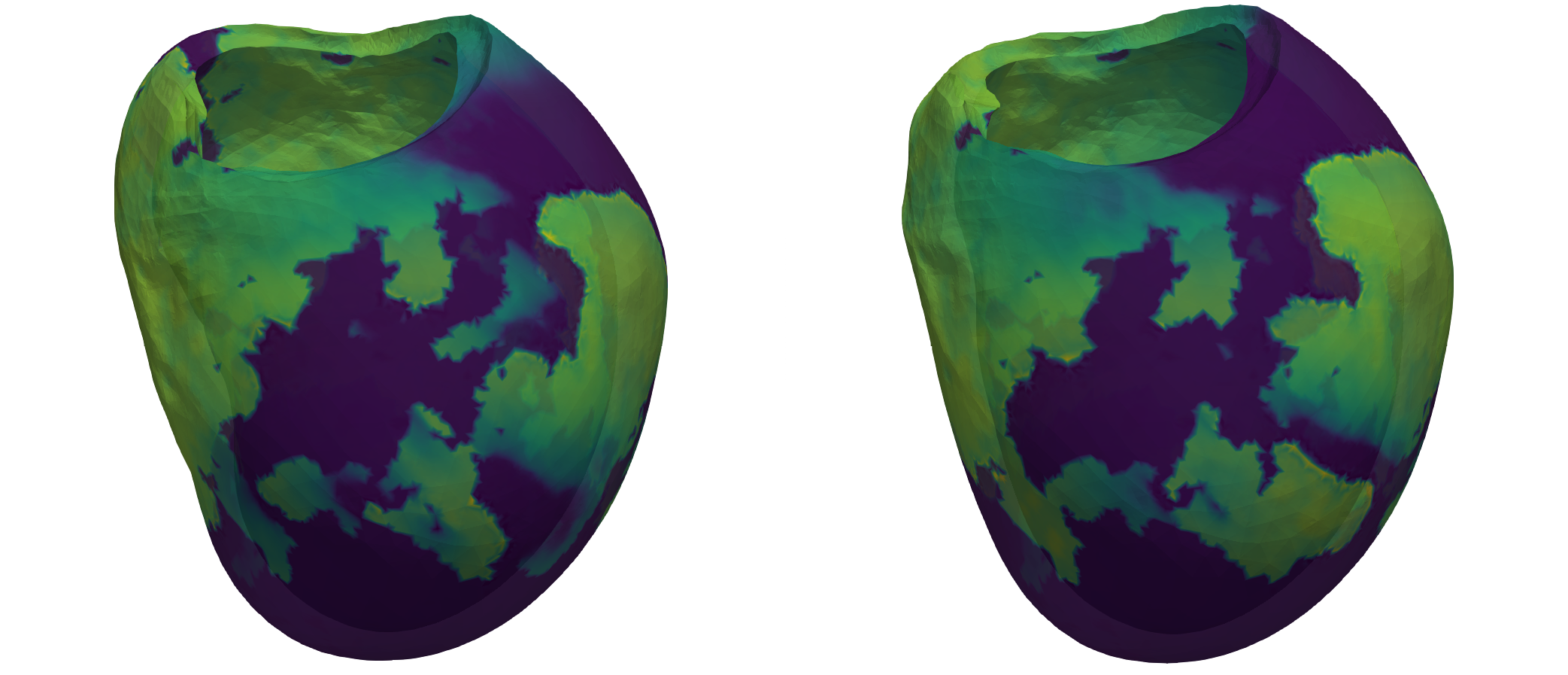}{6.0} \\
	\EMsnapshotVTVT{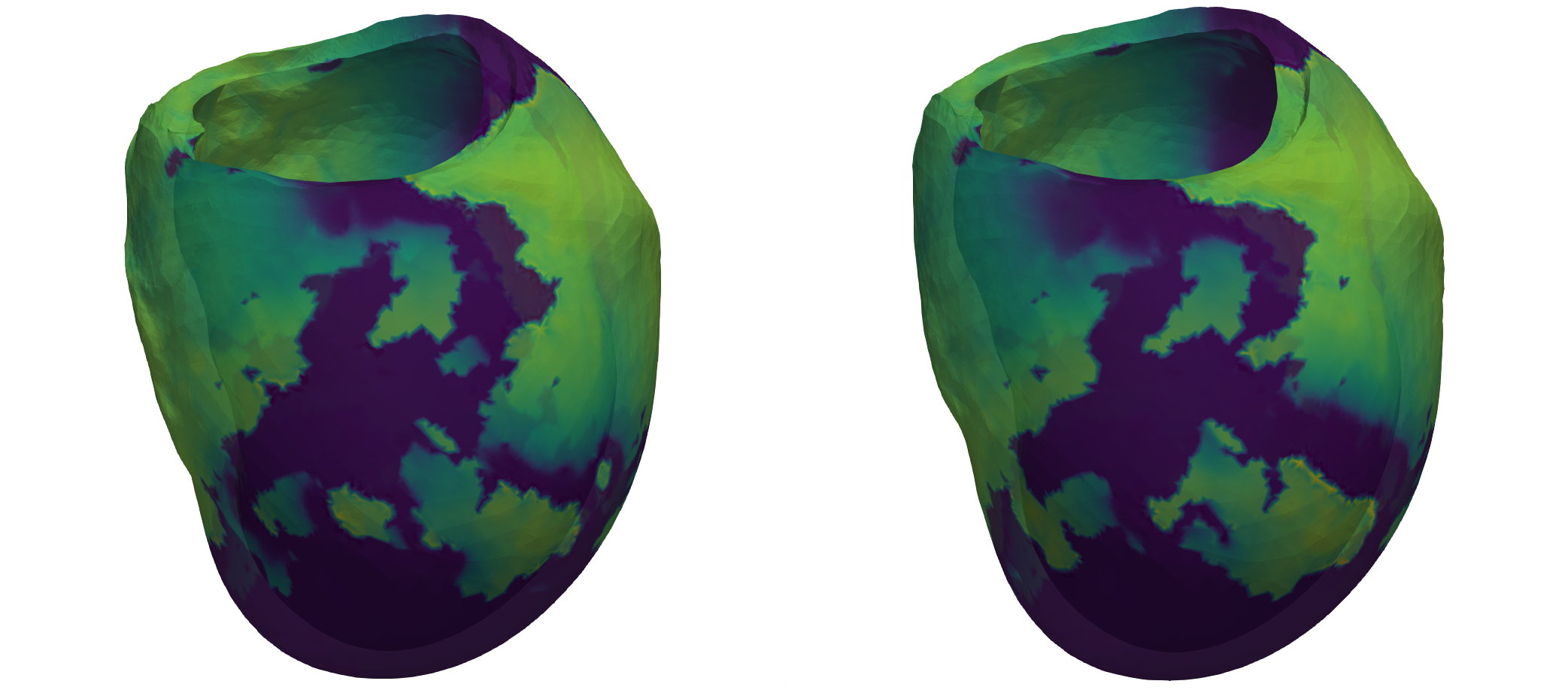}{6.5} \quad
	\EMsnapshotVTVT{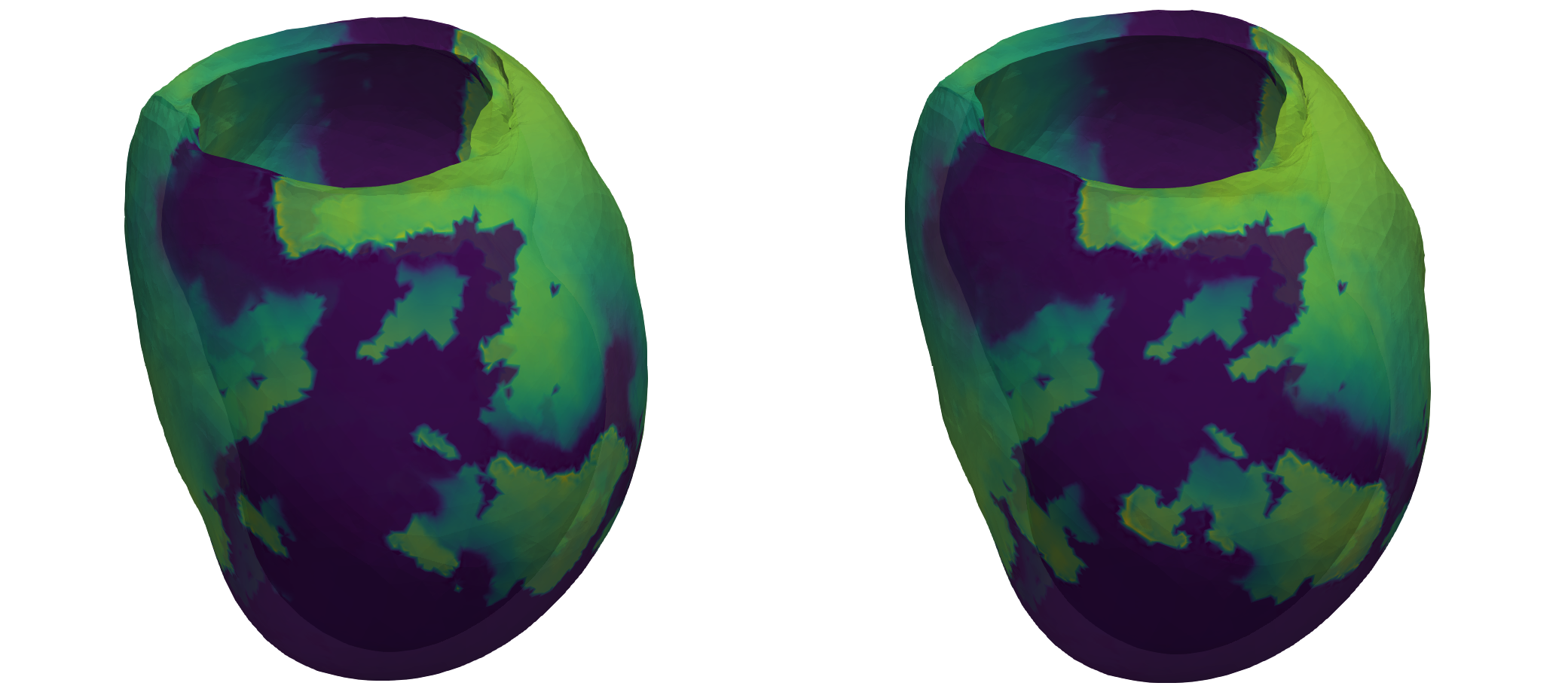}{7.0} \\
	\EMsnapshotVTVT{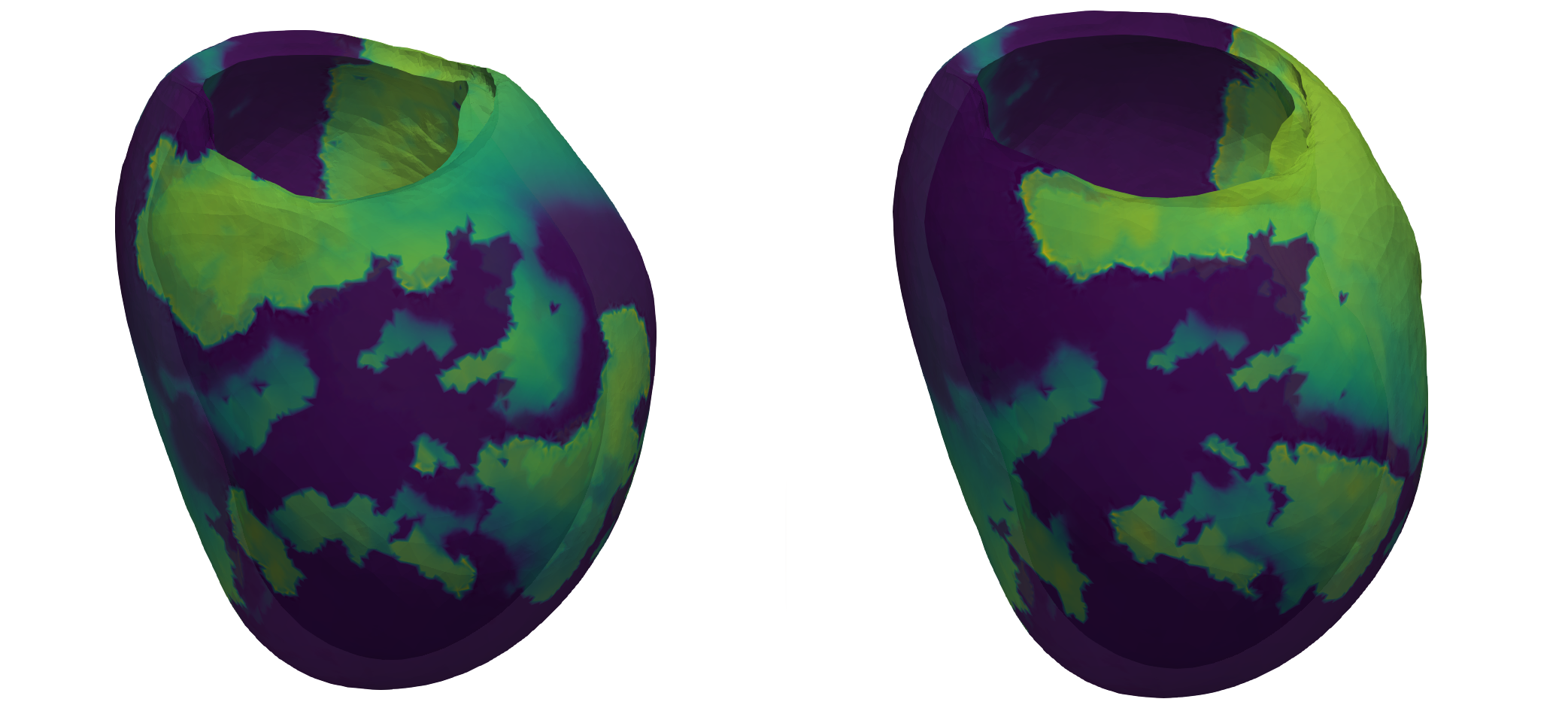}{7.5} \quad
	\EMsnapshotVTVT{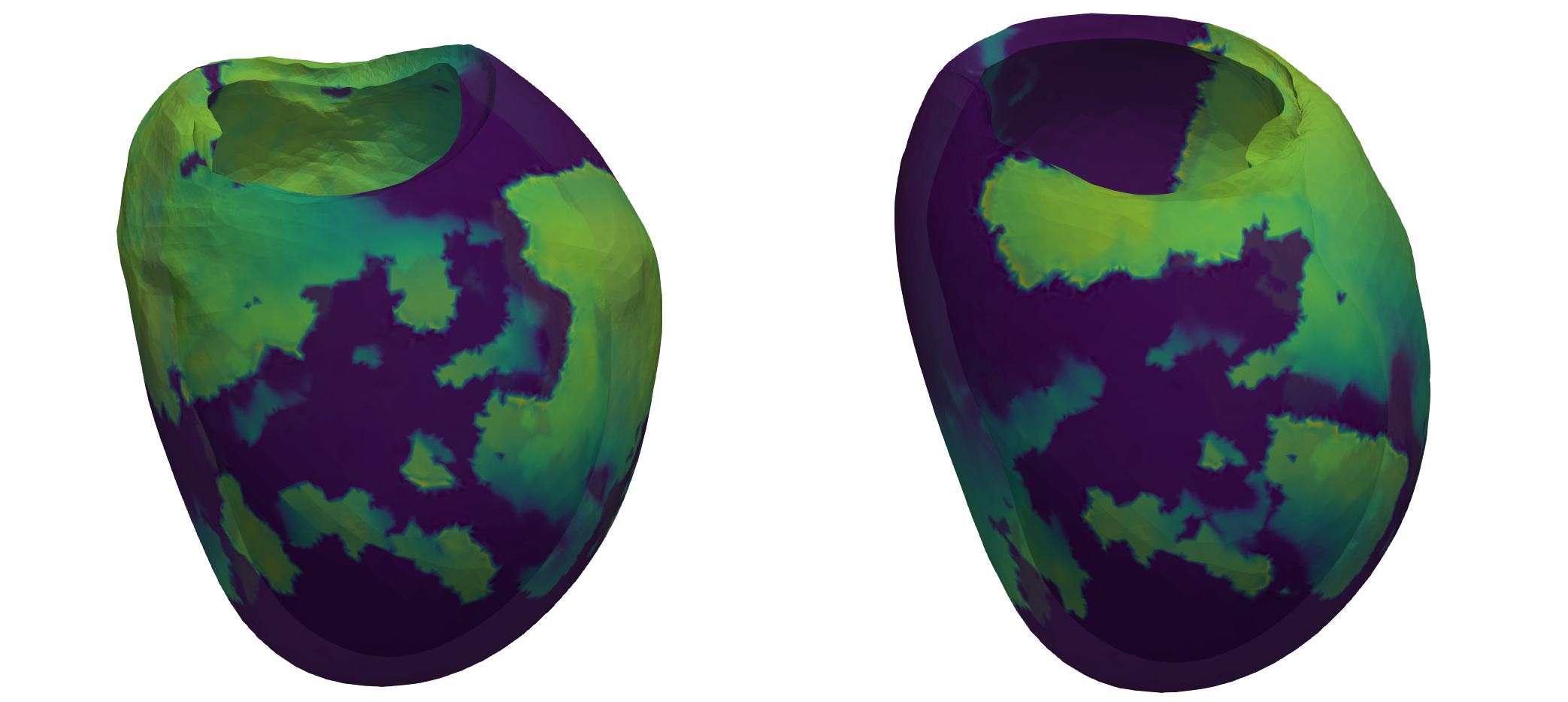}{8.0} \\
	\EMsnapshotVTVT{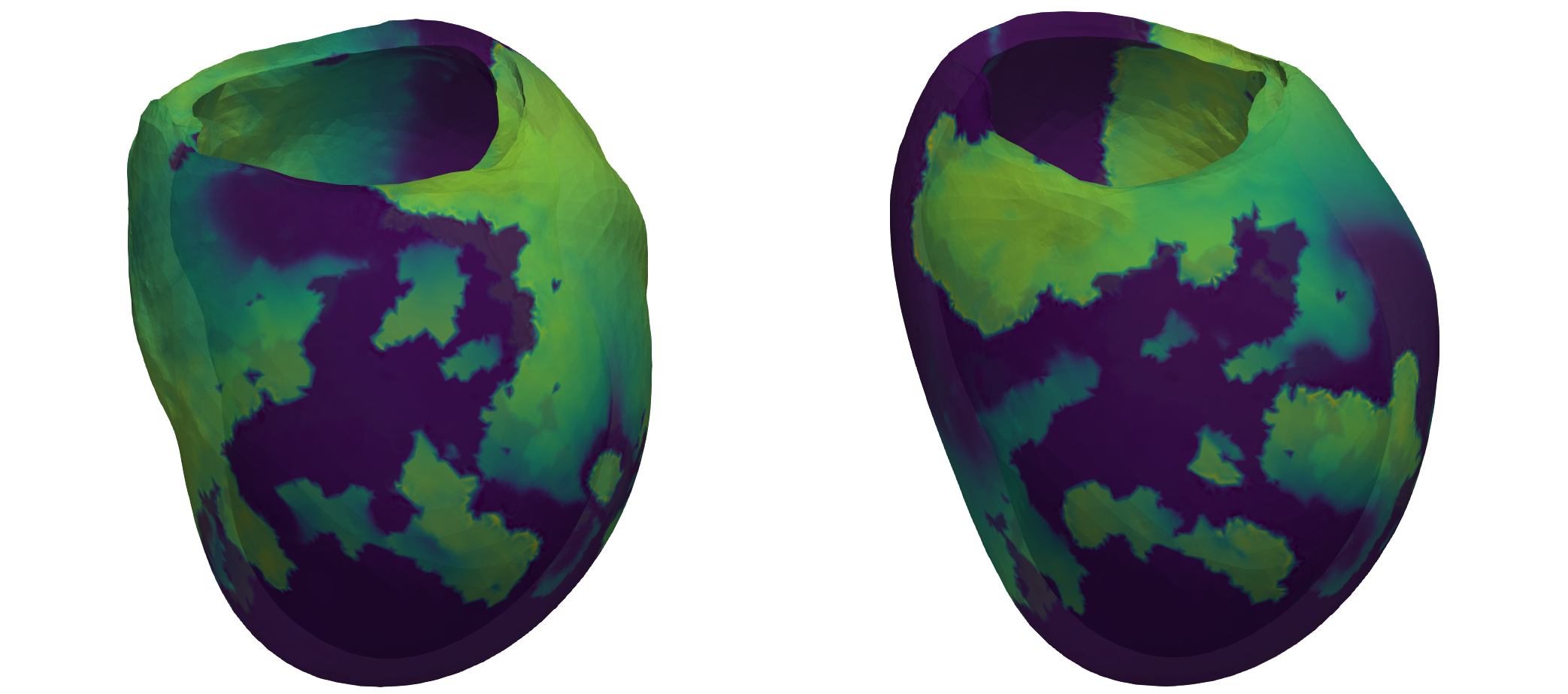}{8.5} \quad
	\EMsnapshotVTVT{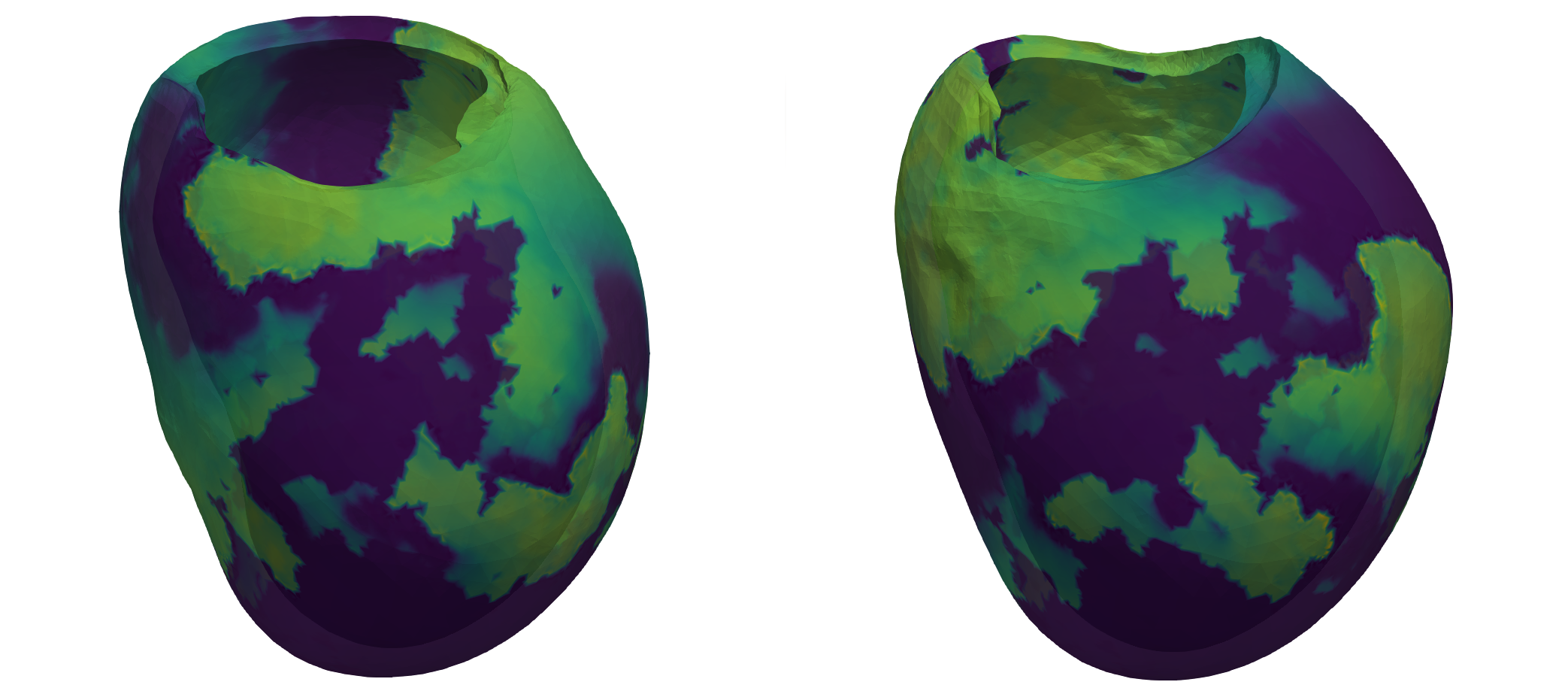}{9.0} \\
	\EMsnapshotVTVT{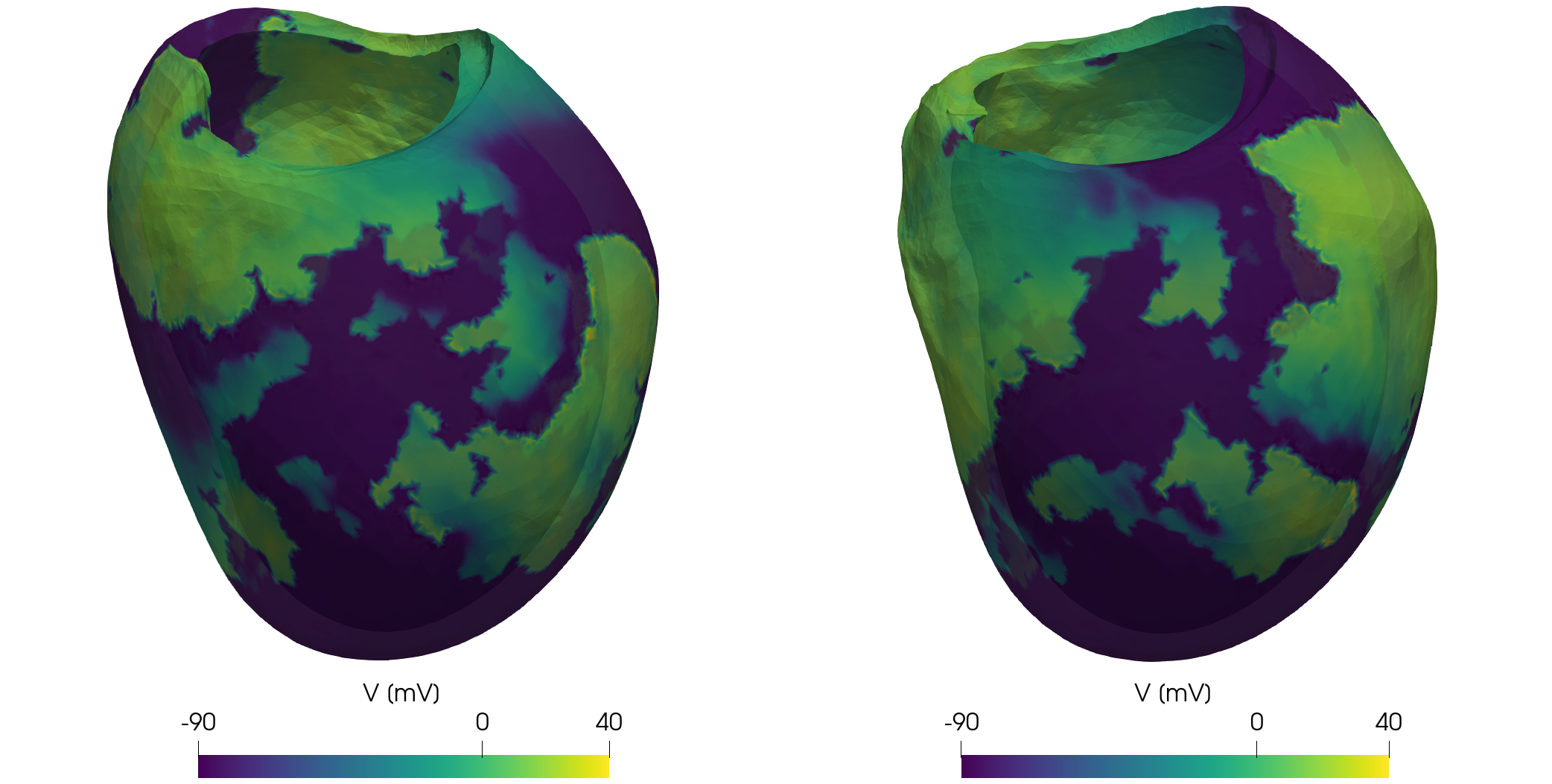}{9.5} \quad
	\EMsnapshotVTVT{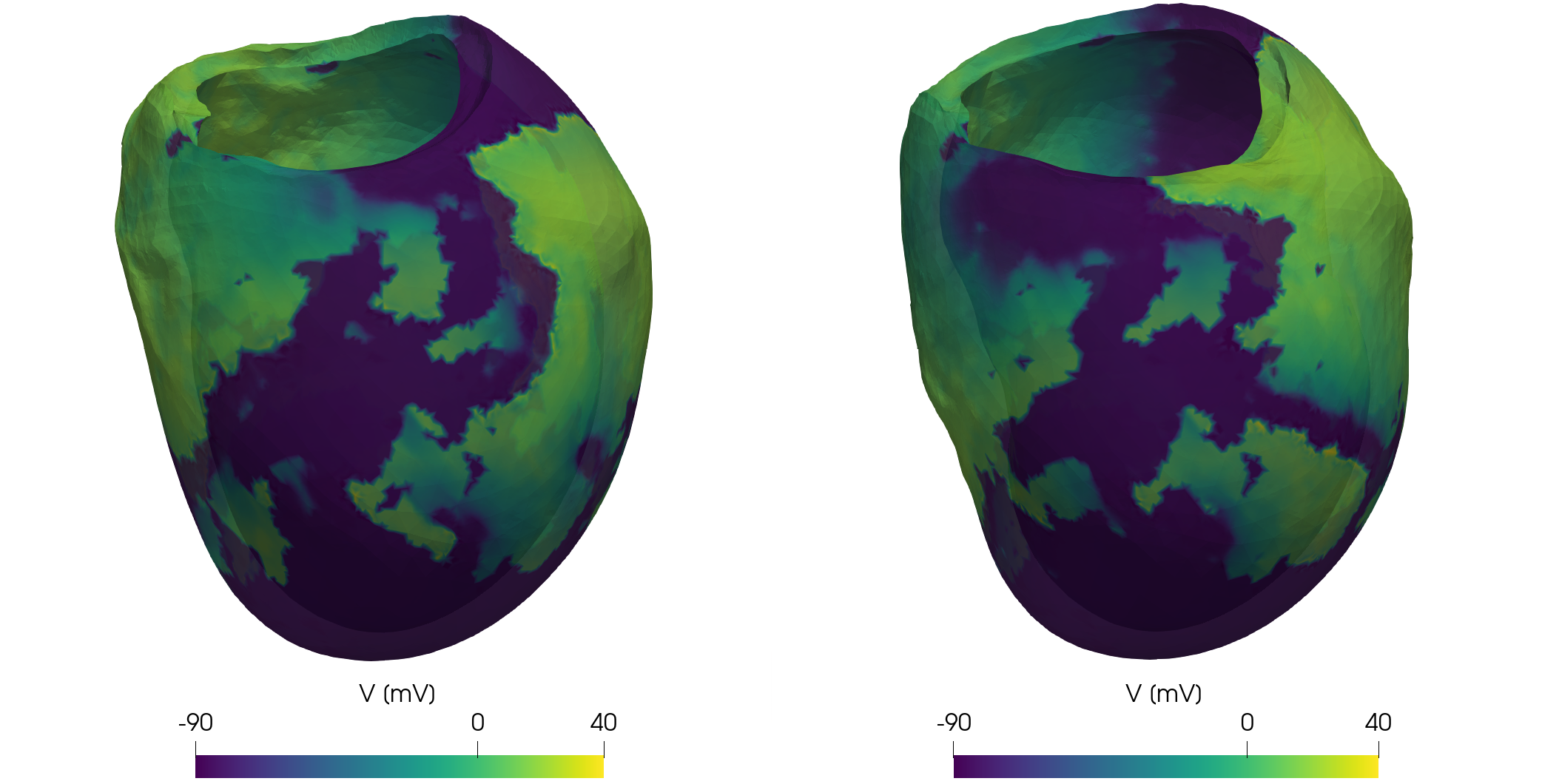}{10.0}

	\caption{Propagation of the transmembrane potential $V=85.7 u - 84$ during VT for a patient-specific LV with ICM. We depict two electromechanics simulations with a healthy parametrization of the circulation model (left) and a pathological one (right). The geometry is warped by the displacement vector in both cases.}
	\label{fig: EM_baseline_simulation_V_displacement_2}
\end{figure}

\begin{figure}[t!]
\advance\leftskip-3.0cm
\includegraphics[keepaspectratio, width=1.5\textwidth]{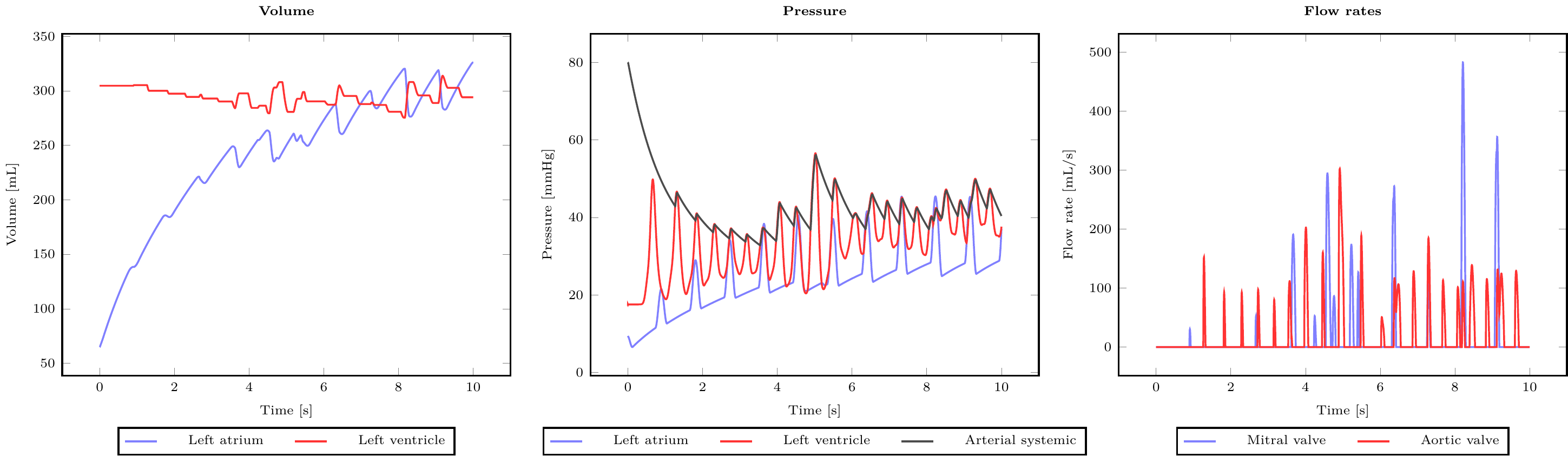}
\caption{Electromechanics simulation of VT for a patient-specific LV with ICM and a healthy parametrization of the circulation model: blood pool volume over time for the 0D left atrium and the 3D LV (left), pressure over time for the left atrium, the left ventricle and the arterial systemic part of the cardiovascular system (center), flow rates of the mitral valve and the aortic valve (right).}
\label{fig: pt339_circulation_VT_healthy}
\end{figure}

\begin{figure}[t!]
\advance\leftskip-3.0cm
\includegraphics[keepaspectratio, width=1.5\textwidth]{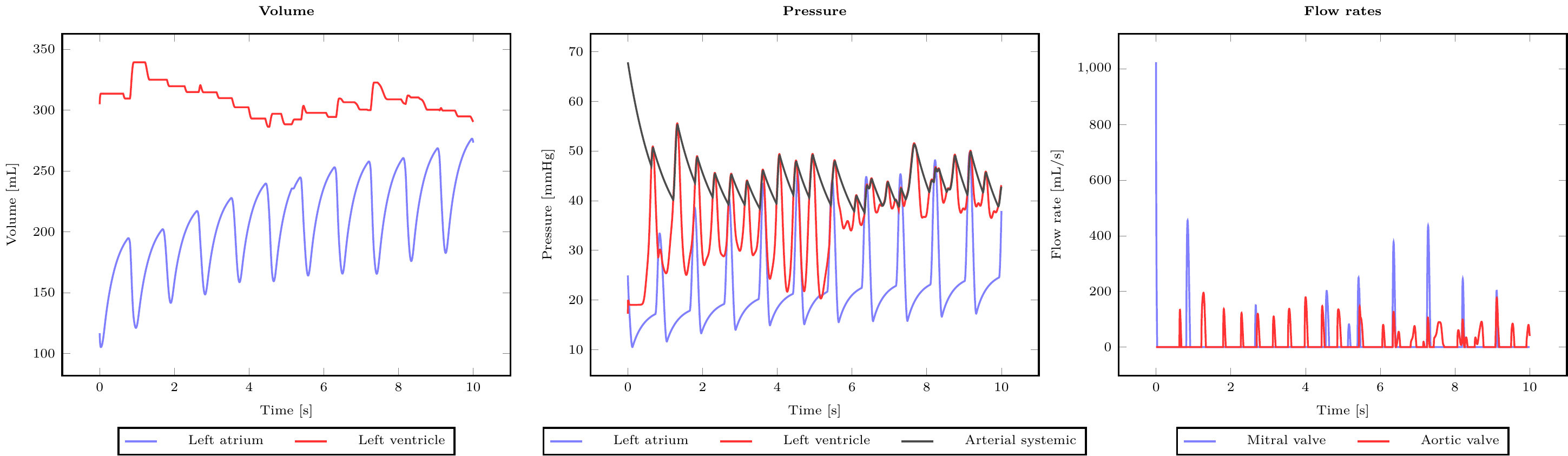}
\caption{Electromechanics simulation of VT for a patient-specific LV with ICM and a pathological parametrization of the circulation model: blood pool volume over time for the 0D left atrium and the 3D LV (left), pressure over time for the left atrium, the left ventricle and the arterial systemic part of the cardiovascular system (center), flow rates of the mitral valve and the aortic valve (right).}
\label{fig: pt339_circulation_VT_pathological}
\end{figure}

We induced a VT by repeatedly stimulating the LV in one specific location. Again, we considered applied currents with intensities $\EPIappReduced(\boldsymbol{x}, t)$, cubic distribution in space, peak $\EPIappReducedMax$ and duration equal to $\EPIappDuration$. First, we run an electrophysiology simulation with the software package CARP (CardioSolv LLC) \cite{Vigmond2003, Vigmond2008} on the original tetrahedral mesh ($h_\text{mean} = 0.35$ mm) and we induced a persistent VT by following the strategy proposed in \cite{Arevalo}. Then, we used our hexahedral mesh ($h_\text{mean} = 1.5$ mm) and we calibrated the conductivities of the monodomain equation to reproduce a similar activation map, as shown in Fig~\ref{fig: pt339_AT}. Next, we induced a persistent VT using a shorter stimulation protocol. In particular, we delivered two S1 stimuli at times $\SI{0}{\second}$ and $\SI{0.450}{\second}$. Then, we delivered a S2 stimulus at time $\SI{0.760}{\second}$ and a S3 stimulus at $\SI{1.040}{\second}$. We let the VT evolve over multiple heartbeats. While the VT observed in CARP has a basis cycle length $BCL \approx \SI{0.390}{\second}$, the one resulting from \texttt{life\textsuperscript{x}} shows a $BCL \approx \SI{0.420}{\second}$.
Then, we run an electromechanics simulation by using the multiphysics model parameters reported in Appendix~\ref{app:params}. More specifically, for the 0D circulation model, we consider the pathological parametrization reported in Tab.~\ref{tab:parameters circulation pathological}, with steady state initial conditions obtained after running the 0D model in SR for 100 heartbeats.
In Fig.~\ref{fig: EM_baseline_simulation_V_displacement_1} we compare the electrophysiology simulation with the electromechanics one, by depicting the distribution of the transmembrane potential over time.
In the time interval before the VT is triggered, no major differences between the two simulations are observed. On the other hand, after the VT is induced, its morphology around the isthmus significantly changes, as well as the conduction velocity of the electric signal. Indeed, the basis cycle length slightly increases, going to $BCL \approx 0.440$ $s$.
Moreover, given a certain stimulation point and a certain sequence of S1-S2-S3 stimuli, we also conclude that in this case, if a VT is observed by considering electrophysiology only, a VT is also inducible by considering electromechanics.
We also compared two electromechanics simulations which are performed by using a very different parametrization for the 0D circulation model: the first one, which is reported in Tab.~\ref{tab:parameters circulation healthy}, resembles the activity of a heart without infarction, whereas the second one, reported in Tab.~\ref{tab:parameters circulation pathological}, defines a pathological cardiovascular system with an infarcted heart. For further details about the second setting, we refer to Appendix~\ref{app:params}.
In Fig.~\ref{fig: EM_baseline_simulation_V_displacement_2} we highlight that different parametrizations of the circulation model induce differents activation patterns and consequently different displacement fields. Nevertheless, we employ the very same sequence S1-S2-S3 stimuli to induce the VT.
In Figs.~\ref{fig: pt339_circulation_VT_healthy} and~\ref{fig: pt339_circulation_VT_pathological} we illustrate the behavior of left atrium pressure and volume, left ventricle pressure and volume, arterial systemic pressure, mitral valve and aortic valve flow rates, during VT. We highlight that pressures in the aorta and in the LV are either dropping or oscillating on unsustainable levels for the cardiovascular system. Indeed, the LV is not able to push blood into the circulation system and the blood pool volume presents small variations over time. On the other hand, both left atrium pressure and volume are increasing during VT. In particular, due to electric isolation between atria and ventricles, the left atrium is assumed to follow the SR pacing during the simulated VT, which lasts for $\SI{10}{\second}$. Moreover, the flow rates of both mitral valve and aortic valve indicate that, even if no regurgitation occurs, there is no proper synchronization between the different cardiac chambers. For this reason, the left heart function is highly compromised. Finally, the PV loop of the LV is not stabilizing over a certain limit cycle that would not cause its impairment. For this reason, we can classify the VT as hemodynamically unstable regardless the parametrization of the 0D circulation model. This VT is critical for the patient and it may lead to SCD.

\section{Discussion}
\label{sec: conclusions}

We presented an individualized computational model of the electromechanical activity in the LV of a patient with ICM, both in SR and VT. We personalized model parameters by means of numerical simulations in SR to fit the available clinical data. Then, we successfully induced a persistent VT and we studied its effects by combining electrophysiological, mechanical and hemodynamical observations. To the best of our knowledge, this is the first time in which a VT has been analyzed by electromechanics numerical simulations in a patient-specific ventricle with ICM.

Our mathematical parametrization incorporates the heterogeneous distribution of scars, grey zones and non-remodeled regions of human ventricles. Different from previous works in literature \cite{Arevalo, Prakosa}, we model both electrophysiological and mechanical properties. While prior state-of-the-art electromechanical numerical simulations seek to model both normal function and pathological conditions, including heart failure \cite{Strocchi2020}, or the impact of drugs \cite{Margara2021}, only the SR case is addressed and tissue heterogeneity of the myocardium is not kept into account over the entire electromechanics pipeline. In our electromechanical model, the coupling with a closed-loop system, as proposed in \cite{Regazzoni2020PartI}, and the numerical scheme developed in \cite{Regazzoni2020PartII}, allow for the effective numerical simulation of VTs. Indeed, this approach does not discriminate among the four different phases of the PV loop as in prior formulations of electromechanics models \cite{Chapelle, Gerbi, Nobile}.

For SR simulations, we observed major differences in the hemodynamics of an LV with ICM with respect to a reference healthy case. Specifically, when incorporating the pathological remodeling, we noticed an increase in EDV and EDP, and a significant reduction in SV/EF and contractility.
Regarding VT simulations, we showed that, if a certain sequence of stimuli delivered in a specific stimulation point induces a VT in electrophysiology simulations, a VT can be also observed with the same settings in electromechanics simulations. On the other hand, geometric MEFs alter the morphology of the VT, along with the overall conduction velocity, which slightly decreases with respect to the electrophysiology simulations. Indeed, the wave propagation is partially influenced by the displacement of the myocardium while running electromechanics simulations. This is due to the presence of the ventricular deformation inside the formulation of the monodomain equation. Geometric and physiological MEFs, such as the recruitment of SACs, are known to have this effect on the conduction velocity restitution curves \cite{Hu2013}. Moreover, geometric MEFs do not generally affect wave stability, while physiological MEFs may determine wavebreaks and the onset of fibrillation \cite{ColliFranzone2017, Panfilov2010}.

With the 0D circulation model, we can also compute the evolution in time of left ventricle pressure and volume, left atrium pressure and volume, arterial system pressure, mitral valve and aortic valve flow rates. With this information, we could classify a VT as either hemodynamically stable or unstable: the former would let the LV stabilize on a PV loop that does not compromise its function, while the latter entails an unstable behavior of the arterial pressure over time, small variations of the blood pool volume of the LV over time, valves flow rates that are not synchronized \cite{Yamada2012}. We conclude that this specific VT is unstable. Having the ability to non-invasively assess a hemodynamically unstable VT is very useful from the clinical perspective. Our electromechanical model allows a thorough understanding of this VT using numerical simulations, whereas investigating this intraprocedurally could be difficult, since the patient might not support hemodynamically unstable VTs.

For this specific patient, very different parametrizations of the 0D closed-loop circulation model did not change the hemodynamical nature of the VT, which remained unstable. This may suggest that the hemodynamics of the patient is linked to the electromechanical substrate. Furthermore, this may have strong clinical implications, in particular when the parameters of the cardiovascular system of the patient are either not known or very uncertain.
The electromechanics simulations presented in this paper can be employed for precision medicine and to gain a deeper understanding of VTs mechanisms thanks to detailed electric, mechanical and hemodynamical descriptions.

\section{Limitations}

First, our current model only considers the effect of geometric MEFs. The effects of physiological MEFs, such as SACs \cite{Hu2013, Kohl2004}, and more accurate geometric MEFs \cite{LevreroFlorencio2020} would need to be addressed as well, to assess whether the conclusions drawn in this paper are still valid. Nevertheless, the impact of the modeled geometric MEFs in the electromechanics simulations entails significant changes with respect to electrophysiology-only simulations. The second limitation of our current study is that only one in silico electrophysiology VT induced from one pacing location is used as the reference. It would be also interesting to understand if there are VTs that can be triggered only by means of electromechanics simulations, and not in electrophysiology ones, and vice-versa. This would require more numerical simulations with different pacing locations distributed all over the ventricle. The third limitation is that the same mesh is used for electrophysiology and mechanics. From the numerical perspective, this entails high computational costs. For this reason, we have to develop new techniques to properly treat heterogeneous tissue on different space resolutions. We believe that this approach would lead to faster numerical simulations without comprimising accuracy, since similar results are observed between different resolutions in electromechanics simulations \cite{Salvador2020}. The fourth limitation is that we have only considered a single LV with ICM instead of a bi-ventricular model. Finally, our approach was evaluated only on one patient. For better generalizability to clinical settings, numerical simulations should be conducted on a larger patient cohort.

\section*{Acknowledgements}
This project has received funding from the European Research Council (ERC) under the European Union's Horizon 2020 research and innovation programme (grant agreement No 740132, iHEART - An Integrated Heart Model for the simulation of the cardiac function, P.I. Prof. A. Quarteroni). M. Salvador acknowledges Dr. S. Pagani and Dr. F. Regazzoni for their advice in cardiac electrophysiology and the useful discussions about cardiac electromechanics.
\begin{center}
	\raisebox{-.5\height}{\includegraphics[width=.15\textwidth]{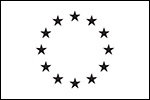}}
	\hspace{2cm}
	\raisebox{-.5\height}{\includegraphics[width=.15\textwidth]{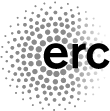}}
\end{center}

\clearpage

\appendix

\section{Model parameters} \label{app:params}
We provide the list of parameters adopted for both SR and VT simulations (Sec.~\ref{sec: sinusrhythmsimulations} and Sec.~\ref{sec: arrhythmiasimulations}). Specifically, Tab.~\ref{tab:parameters electrophysiology} contains the parameters related to the electrophysiology model and Tab.~\ref{tab:parameters mechanics} those related to the mechanics model. For the TTP06 model, we adopt the parameters of the original paper (for myocardial cells) \cite{TTP06}.
Finally, Tab.~\ref{tab:parameters circulation healthy} and Tab.~\ref{tab:parameters circulation pathological} contains the parameters of the circulation model in healthy and pathological conditions, respectively. In particular, for the pathological case we increased the resistance of the arterial system while decreasing its capacitance, to maintain their product almost constant. Morever, we raised the passive elastance of the right ventricle while reducing the active one. Indeed, both from LGE-MRI and Cine MRI, we noticed that the right ventricle of this patient might be affected by ischemic cardiomyopathy as well. Moreover, it presents a small dimension and a low SV ($\approx \SI{25}{\milli\liter}$). Finally, we considered a fibrosis-free dilated left atrium, with higher active elastance and higher resting volume.

\begin{table}[h!]
	\centering
	\begin{tabular}{lrr|lrr}
		\toprule
		Variable & Value & Unit & Variable & Value & Unit \\
		\midrule
		\multicolumn{3}{l|}{\textbf{Conductivity tensor}} & \multicolumn{3}{l}{\textbf{Applied current}} \\
		$\sigma_{\text{l}}$ & \num{0.6714e-4} & \si{\meter\squared\per\second}   &
		$\EPIappReducedMax$  & \num{35}       & \si{\volt\per\second}   \\
		$\sigma_{\text{t}}$ & \num{0.0746e-4}    & \si{\meter\squared\per\second}   &
		$\EPIappDuration$    & \num{5e-3} & \si{\second}   \\
		$\sigma_{\text{n}}$ & \num{0.0746e-4} & \si{\meter\squared\per\second}   &
		& &  \\
		\bottomrule
	\end{tabular}
	\caption{Parameters of the electrophysiology model.}
	\label{tab:parameters electrophysiology}
\end{table}

\begin{table}[h!]
	\centering
	\begin{tabular}{lrr|lrr}
		\toprule
		Variable & Value & Unit & Variable & Value & Unit \\
		\midrule
		\multicolumn{3}{l|}{\textbf{Constitutive law}} & \multicolumn{3}{l}{\textbf{Boundary conditions}} \\
		$B$   & \num{5e4}           & \si{\pascal}   &
		$\BCmecKepiN$      & \num{2e5}       & \si{\pascal\per\meter}   \\
		$C$    & \num{0.88e3}        & \si{\pascal}   &
		$\BCmecKepiT$      & \num{2e5}       & \si{\pascal\per\meter}   \\
		$\bff$    & 8          & $-$   &
		$\BCmecCepiN$      & \num{2e4}       & \si{\pascal\second\per\meter}   \\
		$\bss$    & 6        & $-$   &
		$\BCmecCepiT$   &  \num{2e3}      & \si{\pascal\second\per\meter}   \\
		$\bnn$    & 3          & $-$   &
		& &    \\
		$\bfs$    & 12         & $-$   &
		\multicolumn{3}{l}{\textbf{Activation}} \\
		$\bfn$    & 3         & $-$   &
		$\hat{\mu}_A^1$   &  1.5      & \si{\second \, \micro M^2}     \\
		$\bsn$    & 3         & $-$   &
		$\hat{\mu}_A^2$ &  3       & \si{\second \, \micro M^2} \\
		$\rho_\text{s}$ & $10^3$         & \si{\kilogram \per \cubic\meter}  &
		$\hat{\mu}_A^3$ &  1.2       & \si{\second \, \micro M^2} \\
	    & & &
	    $\hat{\mu}_A^4$ &  5       & \si{\second \, \micro M^2} \\
		& &   \\
		\multicolumn{3}{l|}{\textbf{Windkessel}} \\
	    $C$ &  4.0e-10 & \si{\cubic\meter \per \pascal} \\
   	    $R$ &  5.0e7 & \si{\pascal \, \second \, m^{-3}} \\
		\bottomrule
	\end{tabular}
	\caption{Parameters for activation, mechanics and Windkessel model. For all the other parameters of the active strain model we refer to \cite{Salvador2020}.}
	\label{tab:parameters mechanics}
\end{table}

\begin{table}[h!]
	\centering
	\begin{tabular}{lrr|lrr}
		\toprule
		Variable & Value & Unit & Variable & Value & Unit \\
		\midrule
		\multicolumn{3}{l|}{\textbf{External circulation}} & \multicolumn{3}{l}{\textbf{Cardiac chambers}} \\
		$\RarSYS$   & 0.64           & \si{\mmHg \second \per \milli\liter}   &
		$\EpLA$      & 0.09       & \si{\mmHg \per \milli\liter}   \\
		$\RarPUL$    & 0.032116   & \si{\mmHg \second \per \milli\liter}   &
		$\EpRA$      & 0.07       & \si{\mmHg \per \milli\liter}   \\
		$\RvnSYS$    & 0.035684   & \si{\mmHg \second \per \milli\liter}   &
		$\EpRV$      & 0.05       & \si{\mmHg \per \milli\liter}   \\
		$\RvnPUL$    & 0.1625        & \si{\mmHg \second \per \milli\liter}   &
		$\EaMaxLA$   & 0.07       & \si{\mmHg \per \milli\liter}   \\
		$\CarSYS$    & 1.2          & \si{\milli\liter \per \mmHg}   &
		$\EaMaxRA$   & 0.06       & \si{\mmHg \per \milli\liter}   \\
		$\CarPUL$    & 10.0         & \si{\milli\liter \per \mmHg}   &
		$\EaMaxRV$   & 0.55       & \si{\mmHg \per \milli\liter}   \\
		$\CvnSYS$    & 60.0         & \si{\milli\liter \per \mmHg}   &
		$\VnLA$      & 4.0          & \si{\milli\liter}   \\
		$\CvnPUL$    & 16.0         & \si{\milli\liter \per \mmHg}   &
		$\VnRA$      & 4.0          & \si{\milli\liter}   \\
		$\LarSYS$    & \num{5e-3}   & \si{\mmHg \second\squared \per \milli\liter}  &
		$\VnRV$      & 10.0         & \si{\milli\liter}   \\
		$\LarPUL$    & \num{5e-4}         & \si{\mmHg \second\squared \per \milli\liter}  &
		\multicolumn{3}{l}{\textbf{Cardiac valves}} \\
		$\LvnSYS$    & \num{5e-4}         & \si{\mmHg \second\squared \per \milli\liter}  &
		$\Rmin$      & 0.0075       & \si{\mmHg \second \per \milli\liter}   \\
		$\LvnPUL$    & \num{5e-4}         & \si{\mmHg \second\squared \per \milli\liter}  &
		$\Rmax$      & 75006.2      & \si{\mmHg \second \per \milli\liter}   \\
		\bottomrule
	\end{tabular}
	\caption{Parameters of the circulation model in healthy conditions (mainly taken from \cite{Regazzoni2020PartI, Regazzoni2020PartII}). We always consider a heartbeat period $T = \SI{0.92}{\second}$.}
	\label{tab:parameters circulation healthy}
\end{table}

\begin{table}[h!]
	\centering
	\begin{tabular}{lrr|lrr}
		\toprule
		Variable & Value & Unit & Variable & Value & Unit \\
		\midrule
		\multicolumn{3}{l|}{\textbf{External circulation}} & \multicolumn{3}{l}{\textbf{Cardiac chambers}} \\
		$\RarSYS$   & 1.0           & \si{\mmHg \second \per \milli\liter}   &
		$\EpLA$      & 0.09       & \si{\mmHg \per \milli\liter}   \\
		$\RarPUL$    & 0.032116   & \si{\mmHg \second \per \milli\liter}   &
		$\EpRA$      & 0.07       & \si{\mmHg \per \milli\liter}   \\
		$\RvnSYS$    & 0.26          & \si{\mmHg \second \per \milli\liter}   &
		$\EpRV$      & 0.3       & \si{\mmHg \per \milli\liter}   \\
		$\RvnPUL$    & 0.035684  & \si{\mmHg \second \per \milli\liter}   &
		$\EaMaxLA$   & 0.14      & \si{\mmHg \per \milli\liter}   \\
		$\CarSYS$    & 0.8          & \si{\milli\liter \per \mmHg}   &
		$\EaMaxRA$   & 0.06       & \si{\mmHg \per \milli\liter}   \\
		$\CarPUL$    & 10.0         & \si{\milli\liter \per \mmHg}   &
		$\EaMaxRV$   & 0.4       & \si{\mmHg \per \milli\liter}   \\
		$\CvnSYS$    & 60.0         & \si{\milli\liter \per \mmHg}   &
		$\VnLA$      & 5.0          & \si{\milli\liter}   \\
		$\CvnPUL$    & 16.0         & \si{\milli\liter \per \mmHg}   &
		$\VnRA$      & 4.0          & \si{\milli\liter}   \\
		$\LarSYS$    & \num{5e-3}   & \si{\mmHg \second\squared \per \milli\liter}  &
		$\VnRV$      & 10.0         & \si{\milli\liter}   \\
		$\LarPUL$    & \num{5e-4}   & \si{\mmHg \second\squared \per \milli\liter}  &
		\multicolumn{3}{l}{\textbf{Cardiac valves}} \\
		$\LvnSYS$    & \num{5e-4}   & \si{\mmHg \second\squared \per \milli\liter}  &
		$\Rmin$      & 0.0075       & \si{\mmHg \second \per \milli\liter}   \\
		$\LvnPUL$    & \num{5e-4}         & \si{\mmHg \second\squared \per \milli\liter}  &
		$\Rmax$      & 75006.2      & \si{\mmHg \second \per \milli\liter}   \\
		\bottomrule
	\end{tabular}
	\caption{Parameters of the circulation model in pathological conditions. We always consider a heartbeat period $T = \SI{0.92}{\second}$.}
	\label{tab:parameters circulation pathological}
\end{table}

\clearpage
\bibliographystyle{plain}
\bibliography{references}

\end{document}